\documentclass[11pt,reqno]{article}
\usepackage{amssymb, amsmath, amsthm, amsfonts, amscd, epsfig}
\usepackage[english]{babel}
\usepackage{bm}
\usepackage{graphicx, epsfig}
\newtheorem{thm}{Theorem}[section]
\newtheorem{cor}[thm]{Corollary}
\newtheorem{lem}[thm]{Lemma}
\newtheorem{prop}[thm]{Proposition}

\newcommand{\nm}{\noalign{\smallskip}}

\def\ep{\epsilon}

\newcommand{\Bx}{\mathbf{x}}

\newcommand{\RR}{\mathbb{R}}

\newcommand{\Scal}{\mathcal{S}}

\newcommand{\Kcal}{\mathcal{K}}
\newcommand{\p}{\partial}

\newcommand{\pd}[2]{\frac {\p #1}{\p #2}}
\newcommand{\ds}{\displaystyle}
\newcommand{\eqnref}[1]{(\ref {#1})}

\newcommand{\beq}{\begin{equation}}
\newcommand{\eeq}{\end{equation}}

\newcommand{\la}{\langle}
\newcommand{\ra}{\rangle}

\numberwithin{equation}{section}
\numberwithin{figure}{section}

 \def\p{\partial}
\def \Vh0{\stackrel{\circ}{V}_h}

\def\l{\label}  \def\f{\frac}  

\def\l|{\left|}
\def\r|{\right|}

\newcommand{\lc}
{\mathrel{\raise2pt\hbox{${\mathop<\limits_{\raise1pt\hbox
{\mbox{$\sim$}}}}$}}}

\newcommand{\gc}
{\mathrel{\raise2pt\hbox{${\mathop>\limits_{\raise1pt\hbox{\mbox{$\sim$}}}}$}}}

\newcommand{\ec}
{\mathrel{\raise2pt\hbox{${\mathop=\limits_{\raise1pt\hbox{\mbox{$\sim$}}}}$}}}

\def\be{\begin{equation}} \def\ee{\end{equation}}

\def\bea{\begin{eqnarray}}  \def\eea{\end{eqnarray}}

\def\beas{\begin{eqnarray*}} \def\eeas{\end{eqnarray*}}

\def\bn{\begin{enumerate}} \def\en{\end{enumerate}}

\def\bd{\begin{description}} \def\ed{\end{description}}

 \def\p{\partial} 
\def \Vh0{\stackrel{\circ}{V}_h}

\def\l{\label}  \def\f{\frac}  

\def\p{\partial}

\def\bb{\begin{equation}} \def\ee{\end{equation}}

\def\beqn{\begin{eqnarray}}  \def\eqn{\end{eqnarray}}

\def\beqnx{\begin{eqnarray*}} \def\eqnx{\end{eqnarray*}}

\def\bn{\begin{enumerate}} \def\en{\end{enumerate}}

\def\bd{\begin{description}} \def\ed{\end{description}}

\begin{document}


\title{Shape reconstruction of nanoparticles from their associated  plasmonic resonances\thanks{\footnotesize This work was supported  by the ERC Advanced Grant Project MULTIMOD--267184. Hai Zhang acknowledges a startup fund from HKUST.}} 

\author{
Habib Ammari\thanks{\footnotesize Department of Mathematics, 
ETH Z\"urich, 
R\"amistrasse 101, CH-8092 Z\"urich, Switzerland (habib.ammari@math.ethz.ch, sanghyeon.yu@math.ethz.ch). }
\and  Mihai Putinar\thanks{
Department of Mathematics,  
 University of California at Santa Barbara,  
 Santa Barbara, CA 93106-3080, USA (mputinar@math.ucsb.edu), and
 School of Mathematics \& Statistics, Newcastle University Newcastle upon Tyne, NE1 7RU, United Kingdom (mihai.putinar@ncl.ac.uk)}
 \and  Matias Ruiz\thanks{\footnotesize Department of Mathematics and Applications,
Ecole Normale Sup\'erieure, 45 Rue d'Ulm, 75005 Paris, France
(matias.ruiz@ens.fr).} 
\and Sanghyeon Yu\footnotemark[2]
\and  
Hai Zhang\thanks{\footnotesize 
Department of Mathematics, 
 HKUST,  Clear Water Bay, Kowloon, Hong Kong (haizhang@ust.hk).}
}
\date{}
\maketitle

\begin{abstract}
We prove by means of a couple of examples that plasmonic resonances can be used on one hand to classify shapes of nanoparticles with real algebraic boundaries and on the other hand to reconstruct the separation distance between two nanoparticles from measurements of their first collective plasmonic resonances.  To this end, we explicitly compute the spectral decompositions of the Neumann-Poincar\'{e} operators associated with a class of quadrature domains and two nearly touching disks. Numerical results are included in support of our main findings. 
\end{abstract}

\bigskip

\noindent {\footnotesize Mathematics Subject Classification
(MSC2000): 35R30, 35C20.}


\section{Introduction and main results}

The present paper is a part of an ample and recent effort to 
understand the mathematical structure of inverse problems arising
in nanophotonics. Although very classical, the spectral analysis of the Neumann-Poincar\'e operator emerges as the main theme of investigation.

Consider a domain $\Omega$ with $\mathcal{C}^{1,\eta}$ boundary in $\mathbb{R}^2$ for $\eta>0$. Let $\nu$ denote the outward normal to $\partial \Omega$. Suppose that $\Omega$ contains the origin $0$. For $\alpha=(\alpha_1,\alpha_2) \in \mathbb{N}^2,$ we denote by $\partial_\alpha = \partial_1^{\alpha_1}  \partial_d^{\alpha_2}$ and $\alpha!=\alpha_1! \alpha_2!$. 

The Neumann-Poincar\'{e} operator $\mathcal{K}_{\Omega}^*$ associated with $\Omega$  is defined as follows:
$$
\mathcal{K}_{\Omega}^* [\varphi](x) =\frac{1}{2\pi} \int_{\p \Omega }  \frac{\langle x-y,\nu_x\rangle}{|x-y|^2} \varphi(y) d\sigma(y) ,   \quad x \in \p \Omega.
$$
It is related to the single layer potential $\mathcal{S}_{\Omega}$
given by 
$$
\mathcal{S}_{\Omega} [\varphi](x) =\frac{1}{2\pi} \int_{\p \Omega }  \ln|x-y| \varphi(y) d\sigma(y) ,   \quad x \in \p \Omega,
$$
by the following jump relation:
\begin{align}
\frac{\p\mathcal{S}_\Omega[\varphi]}{\p\nu}\Big|_{\pm} &= (\pm\frac{1}{2} I +\mathcal{K}_\Omega^*)[\varphi] \quad \mbox{for } \varphi \in L^2(\partial \Omega). \label{eqn_jump_single}
\end{align}
 It can be shown the operator $\lambda I - \mathcal{K}_\Omega^*: L^2(\partial \Omega)
\rightarrow L^2(\partial \Omega)$ is invertible for any $|\lambda| > 1/2$. Furthermore, $\mathcal{K}_\Omega^*$ is compact, can be symmetrized in a proper energy space and consequently 
its spectrum is discrete and contained in $]-1/2, 1/2]$; see for instance \cite{book2} for more details.

Using the quasi-static limit of electromagnetic fields, plasmonic resonances are associated with the set of eigenvalues $\lambda_j$ of the Neumann-Poincar\'e operator $\mathcal{K}_\Omega^*$ for which
$\langle \varphi_j, x_i \rangle_{L^2(\partial \Omega)}\neq 0$, for either $i=1$ or $i=2$,  where $\varphi_j$ is an eigenfunction associated to $\lambda_j$ \cite{matias,matias2}. We refer the reader to 
\cite{sisc, pierre, matias, matias2, kang1, hyeonbae, triki, Gri12, helsing, jam} for recent and interesting mathematical results on plasmonic resonances for nanoparticles.

In the present paper we prove that based on plasmonic resonances we can on one hand classify the shape of a class  of domains with real algebraic boundaries and on the other hand recover the separation distance between two components of multiple connected domains. These results have important applications in nanophotonics. They can be used in order to identify the shape and separation distance between plasmonic nanoparticles having known material parameters from measured plasmonic resonances, for which the scattering cross-section is maximized \cite{matias,matias2}.

A real algebraic curve is the zero level set of a bivariate polynomial. Domains enclosed by real algebraic curves (henceforth simply called {\it algebraic domains}) are dense, in Hausdorff metric among all planar domains. On a simpler note, every smooth curve can be approximated by a sequence of algebraic curves. This observation turns algebraic curves into an efficient tool for describing shapes \cite{vetterli,i3e1,i3e2}. 
Note that an algebraic domain which is the sub level set of a polynomial of degree $n$ can uniquely be determined from its set of two-dimensional moments of order less than or equal to $3n$ \cite{ip, dcg}. In this paper we consider a class of algebraic curves determined via conformal mappings by two parameters $m$ and $\delta$,  with 
$m$ being the order of the polynomial parametrizing the curve and $\delta$ being a generalized radius, see (\ref{desc}). One can think of algebraic domains as non-generic, but dense, among all planar domains, as much as polynomials are non-generic, but dense
among all continuous functions on a compact set. In either case, the identifications/reconstructions have to be complemented by a fine analysis of the rate of convergence.

The main results of the present paper are: 

(i) Algebraic domains described by (\ref{desc}) have only two plasmonic resonances asymptotically 
(in $\delta$). Based on these two plasmonic resonances, one can classify them; 

(ii) Two nearly touching disks have an infinite number of plasmonic resonances and the separating distance  can be determined from the measurement of the first plasmonic resonance.

The paper is organized as follows. In section \ref{sect2} we first introduce contracted generalized polarization tensors. 
Then we
give explicit calculations of the Neumann-Poincar\'e operator associated with an algebraic domain. Moreover, we analyze its asymptotic behavior as $\delta$ approaches zero. We compute the first- and second-order contracted polarization tensors, and show how to use them to determine the two parameters describing the algebraic boundaries. 
In section \ref{sect3} we consider two nearly touching disks. We use the bipolar coordinates to compute the spectrum of the associated Neumann-Poincar\'e operator. We show that all the eigenvalues of the associated Neumann-Poincar\'e operator contribute to the set of  plasmonic resonances. From the first-order polarization tensor, we show that we can recover the separating distance between the disks. In section \ref{sect4} we illustrate our main findings in this paper with several numerical examples. In  appendix \ref{appendixA}, we estimate the blow up of the gradient of the potential for two nearly touching disks at plasmonic resonances. This result generalizes to the plasmonic case, the estimates derived in \cite{akl, jde} for nearly touching disks with degenerate conductivities.  


%

\section{Plasmonic resonance for algebraic domains} \label{sect2}

\subsection{Contracted generalized polarization tensors}
Given a harmonic function $H$ in the whole plane, we consider the following transmission problem:
\beqn
    \begin{cases}
        \nabla \cdot (\sigma_{\Omega} \nabla u) = 0 &\text{ in }\; \mathbb{R}^2, \\[1.5mm]
         u - H = O(|x|^{-1}) &\text{ as }\; |x| \rightarrow \infty,
    \end{cases}
    \label{transmission}
\eqn
where $\sigma_{\Omega}  = \sigma \chi(\Omega) +  \chi(\mathbb{R}^2 \backslash \overline{\Omega})$ with $\sigma>0$, and $\chi(\Omega)$ and $\chi(\mathbb{R}^2 \backslash \overline{\Omega})$ are  the characteristic functions of $\Omega$ and $\mathbb{R}^2 \backslash \overline{\Omega}$, respectively.  From \cite{book1}, we have
\beqn
    u = H+ \mathcal{S}_{\Omega} (\lambda I - \mathcal{K}_\Omega^*)^{-1} [\frac{\partial H}{\partial \nu}]  \, ,
    \label{scattered}
\eqn
where \be \label{deflambda} 
\lambda= \frac{\sigma+1}{2(\sigma-1)}.\ee  Decomposition (\ref{scattered}) of $u$ together with
\beqn
    H (x) = \sum_{\alpha \in \mathbb{N}^d } \f{1}{\alpha !} \partial^\alpha H(0) x^\alpha
\label{series}
\eqn
 and 
$$
\ds \Gamma(x-y) = \sum_{|\beta|=0}^{+\infty} \frac{(-1)^{|\beta|}}{\beta!} \partial^\beta_x \Gamma(x) y^\beta, \quad y \mbox{ in a compact set}, \, |x| \rightarrow +\infty, 
$$
where $\Gamma$ is the fundamental solution to the Laplacian, 
yields the far-field behavior  \cite[p. 77] {book1}
\beqn
    (u - H)(x) =  \sum_{|\alpha|, |\beta| \geq 1 } \f{1}{\alpha ! \beta!} \partial^\alpha H(0) \bigg[ \int_{\partial \Omega} y^\beta (\lambda I - \mathcal{K}_\Omega^*)^{-1} [\frac{\partial x^\alpha}{\partial \nu}](y)  d \sigma(y) \bigg]\; \partial^\beta \Gamma(x) 
    \label{scattered2}
\eqn
as $|x| \rightarrow + \infty$.
Introduce the {\it generalized polarization tensors} \cite{book1}: 
$$M_{\alpha \beta}(\lambda, \Omega) : =  \int_{\partial \Omega} y^\beta (\lambda I - \mathcal{K}_\Omega^*)^{-1} [\frac{\partial x^\alpha}{\partial \nu}](y)\, d\sigma(y), \qquad \alpha, \beta \in \mathbb{N}^d.$$
We call $M := M_{\alpha\beta}$ for $|\alpha|=|\beta|=1$ the first-order polarization tensor.

For a positive integer
$m$, let $P_m(x)$ be the complex-valued polynomial
\begin{equation}
P_m(x) = (x_1 + ix_2)^m := \sum_{|\alpha| = m} a^m_\alpha x^\alpha
+ i\sum_{|\beta| = m} b^m_\beta x^\beta. \label{eq:Pdef}
\end{equation}
Using polar coordinates $x = re^{i\theta}$, the above coefficients
$a^m_\alpha$ and $b^m_\beta$ can also be characterized by
\begin{equation}
\sum_{|\alpha| = m} a^m_\alpha x^\alpha = r^m \cos m\theta, \text{
and } \sum_{|\beta| = m} b^m_\beta x^\beta = r^m \sin m\theta.
\label{eq:abcomp}
\end{equation}
We introduce the {\it contracted generalized polarization tensors} to be the
following linear combinations of generalized polarization tensors using the coefficients in
\eqref{eq:Pdef}:
\begin{align*}
M^{cc}_{mn} = \sum_{|\alpha| = m} \sum_{|\beta| = n} a^m_\alpha a^n_\beta M_{\alpha \beta}, \quad
M^{cs}_{mn} = \sum_{|\alpha| = m} \sum_{|\beta| = n} a^m_\alpha b^n_\beta M_{\alpha \beta},\\
M^{sc}_{mn} = \sum_{|\alpha| = m} \sum_{|\beta| = n} b^m_\alpha a^n_\beta M_{\alpha \beta},\quad
M^{ss}_{mn} = \sum_{|\alpha| = m} \sum_{|\beta| = n} b^m_\alpha
b^n_\beta M_{\alpha \beta}. 
\end{align*}
It is clear that 
\begin{align*}
M^{cc}_{mn} = \int_{\partial \Omega} \Re \{ P_n\}  (\lambda I - \mathcal{K}_\Omega^*)^{-1} [\frac{\partial \Re \{ P_m\}}{\partial \nu}]\, d\sigma,  \\
M^{cs}_{mn} = \int_{\partial \Omega} \Im \{ P_n\}  (\lambda I - \mathcal{K}_\Omega^*)^{-1} [\frac{\partial \Re \{ P_m\}}{\partial \nu}]\, d\sigma,\\
M^{sc}_{mn} = \int_{\partial \Omega} \Re \{ P_n\}  (\lambda I - \mathcal{K}_\Omega^*)^{-1} [\frac{\partial \Im \{ P_m\}}{\partial \nu}]\, d\sigma,\\
M^{ss}_{mn} = \int_{\partial \Omega} \Im \{ P_n\}  (\lambda I - \mathcal{K}_\Omega^*)^{-1} [\frac{\partial \Im \{ P_m\}}{\partial \nu}]\, d\sigma. 
\end{align*}
We refer to \cite{book2} for further details

As recently shown  \cite{focm, book2, numer}, the contracted generalized polarization tensors can efficiently be used for domain classification. They provide a natural tool for  describing shapes. 
In imaging applications, they can be stably reconstructed from the data by solving a least-squares problem. They capture high-frequency shape oscillations as well as topology. High-frequency oscillations of the shape of a domain are only contained in its high-order contracted generalized polarization tensors.  

The contracted generalized polarization tensors also satisfy simple invariance properties under translation, rotation, and scaling \cite{focm, acha}. 
Based on those properties, a dictionary matching algorithm can be developed. Assuming  that the unknown shape of the target is an exact copy of some element from the dictionary, up to a rigid transform and dilatation, one can identify the target in the dictionary from its contracted generalized polarization tensors with a low computational cost \cite{focm}. If the material parameters of the domain are frequency dependent, then measurements taken at multiple frequencies allow very stable recognition of the targets. In \cite{pnas}, a classification approach based on the spectral properties of only the first-order polarization tensor was proposed and successfully implemented. 

\subsection{Algebraic domains of class $\mathcal{Q}$} 
Let $D$ be the unit disk in $\mathbb{C}$.
For $m\in\mathbb{N}$ and $a\in\mathbb{R}$, define $\Phi_{m,a}:\mathbb{C}\setminus \overline{D}\rightarrow\mathbb{C}$  by
$$
\Phi_{m,a}(\zeta) = \zeta + \frac{a}{\zeta^m}.
$$
Assume that $\Phi_{m,a}$ is injective on $\mathbb{C}\setminus\overline{D}$.
We introduce the class $\mathcal{Q}$ as the collection of all bounded domains $\Omega \subset \mathbb{C}$ bounded by the curves
$$
\p \Omega = \{ \Phi_{m,a}(\zeta): \quad |\zeta|=r_0\} \quad
\mbox{for some } r_0>1, \quad m\in \mathbb{N} \mbox{ and } a \in \mathbb{R}.
$$
Note that $\Phi_{m,a}$ is a conformal mapping from $ \{|\zeta|> r_0\}$ onto $\mathbb{C}\setminus \overline{\Omega}$. In what follows, we shall suppress the subscript $m, a$ from $\Phi_{m,a}$ for the ease of notation.  

Conformal images of the unit disc by rational functions are also called {\it quadrature domains}. We refer to \cite{GP-qd,QD} for details and ramifications of the theory of quadrature domains.
In particular, up to the inversion $z \mapsto 1/z$, the complements of the domains in class $\mathcal{Q}$ are quadrature domains.
We write for convenience $\zeta=e^{\rho+i\theta}$. Let $\rho_0$ be such that $r_0=e^{\rho_0}$. 
Let $J$ be the Jacobian defined by
$$
J=\big|\p_\xi\big(\Phi(e^{\xi})\big)|_{\xi=\rho+i\theta}\big|.
$$
In the $(\rho,\theta)$ plane, the normal derivative $\p/\p\nu$ on $\p\Omega$ is represented as
$$
\frac{\p}{\p\nu}=\frac{1}{J}\frac{\p}{\p\rho}.
$$
Moreover, the boundary $\p\Omega$ is parametrized by
$$
\theta \mapsto \Phi(e^{\rho_0+i\theta})=e^{\rho_0+i\theta}+a e^{-m\rho_0-im\theta}.
$$
If we fix the constant $a$ and change $\rho_0$, then the size and the shape of $\p\Omega$ will change accordingly. In order to leave the shape unchanged, we need to represent the constant $a$ in a different way. We write
\beq\label{parameter_delta}
a=e^{(m+1)\rho_0} \delta.
\eeq
Then the boundary $\p\Omega$ can be represented as 
\begin{equation} \label{desc}
\theta \mapsto \Phi(e^{\rho_0+i\theta})=e^{\rho_0}(e^{i\theta}+\delta  e^{-im\theta}).
\end{equation}
Now, if we fix the constant $\delta$ and change $\rho_0$, then it is clear that only the size changes and the shape stays unaffected.  The parameter $e^{\rho_0}$ can be considered as a generalized radius of $\Omega$ because it determines the size. In conclusion, the shape of $\Omega$ is determined by the two parameters $m$ and $\delta$, while the size by the parameter $\rho_0$.

\subsection{Explicit computation of the Neumann-Poincar\'{e} operator}
In this section, we compute the Neumann-Poincar\'{e} operator on $\p\Omega$ explicitly. We need to compute $\mathcal{K}^*_{\Omega}[J^{-1} \cos{n\theta}]$ and $\mathcal{K}^*_{\Omega}[J^{-1} \sin{n\theta}]$ explicitly.
Our strategy is as follows. Let $u=\mathcal{S}_\Omega[J^{-1} \cos{n\theta}]$ and $v=\mathcal{S}_\Omega[J^{-1}\sin n\theta]$. If $u,v$ can be obtained explicitly, then  $\mathcal{K}^*_\Omega[J^{-1} \cos{n\theta}]$ and $\mathcal{K}^*_\Omega[J^{-1} \sin{n\theta}]$ are immediately derived by using the following identity:
\be\label{eqn_NP_single_identity}
\mathcal{K}^*_\Omega[\varphi]=\frac{1}{2} \bigg(\frac{\p\mathcal{S}[\varphi]}{\p\nu}\Big|_+ +\frac{\p\mathcal{S}[\varphi]}{\p\nu}\Big|_-\bigg),
\ee
which follows from \eqref{eqn_jump_single}.
For simplicity, we consider only $u$.
By using the continuity of the single layer potential and the jump relation \eqref{eqn_jump_single}, we can see that the function $u$ is the solution to the following problem:
\begin{equation}
 \ \left \{
 \begin{array} {ll}
\ds \Delta u= 0 \quad &\mbox{ in } \mathbb{C}\setminus\partial {\Omega},
\\[0.5em]
\ds u|_-=u|_+ \quad &\mbox{ on } \p\Omega, 
\\[0.5em]
\ds \frac{\p u}{\p\nu}\Big|_+-\frac{\p u}{\p\nu}\Big|_- = J^{-1} \cos{n \theta}\quad &\mbox{ on } \p\Omega, 
\\[1em]
\ds u=O(|z|^{-1}) \quad &\mbox{ as } |z|\rightarrow \infty.
 \end{array}
 \right.
 \end{equation}
 Let $\widetilde{u}(\rho,\theta)=(u\circ\Phi)(e^{\rho+i\theta})$. Since $\Phi(\zeta)$ is conformal on $|\zeta|>e^{\rho_0}$, the above problem can be rewritten as follows:
\begin{equation}\label{eqn_tilde_u_1}
 \ \left \{
 \begin{array} {ll}
\ds \Delta u= 0 \quad &\mbox{ in } \mbox{ for } \rho < \rho_0,
\\[0.5em]
\ds \Delta \widetilde u= 0 \quad &\mbox{ for } \rho> \rho_0,
\\[0.5em]
\ds \widetilde{u}|_-=\widetilde{u}|_+ \quad &\mbox{ on } \rho=\rho_0, 
\\[0.5em]
\ds \frac{\p \widetilde{u}}{\p\rho}\Big|_+-\frac{\p \widetilde{u}}{\p\rho}\Big|_- =  \cos{n \theta}\quad &\mbox{ on } \rho=\rho_0, 
\\[1em]
\ds \widetilde u=O(e^{-\rho}) \quad &\mbox{ as } \rho\rightarrow \infty.
 \end{array}
 \right.
 \end{equation}
Note that in (\ref{eqn_tilde_u_1}), the first equation for $u|_{\Omega}$ is  not represented in terms of $\widetilde{u}$. This is due to the singularity of $\Phi(\zeta)$ near $\zeta=0$. Hence, we need to consider $u|_{\Omega}$ more carefully. 
If $a=1$ and $m=1$, then $\Omega$ becomes an ellipse and $(\rho,\theta)$ are called the elliptic coordinates. In this case,  equation \eqref{eqn_tilde_u_1} for $\widetilde{u}$ can be easily solved by imposing some appropriate conditions on $\rho= \rho_0$ and $\rho=0$. 
 However, for general shaped domains, this is not easy. 
 
 Fortunately, we can overcome this difficulty by the fact that the shape of the domain $\Omega$ 
 is defined by a rational function $\Phi(\zeta)=\zeta+a/\zeta^m$. Our strategy is to seek a solution to \eqref{eqn_tilde_u_1} such that 
 $$
 u(z) = \Re\{\mbox{a polynomial of degree $n$ in $z$} \} \quad \mbox{for }z\in\Omega.
 $$
 We can show that, for $1\leq n\leq m$, $u|_\Omega$ in equation \eqref{eqn_tilde_u_1} can be explicitly solved by using the following ansatz:
  \begin{align}
 u|_{\Omega}(z) \propto \Re\{z^n\}&=\Re\Big\{\Big(\zeta+\frac{a}{\zeta^m}\Big)^n \Big\}
 \nonumber\\
 &=\Re\sum_{k=0}^n \begin{pmatrix}n \\ k\end{pmatrix} \zeta^{n-k} \Big(\frac{a}{\zeta^m}\Big)^k
\qquad (\zeta = e^{\rho I \theta}) \nonumber\\
 &=e^{n\rho}\cos n \theta+\sum_{k=1}^n a^k\begin{pmatrix}n \\ k\end{pmatrix}e^{-t^{mn}_k\rho}\cos t^{mn}_k\theta,\label{eqn_u_Omega_zn}
 \end{align}
where the constant $t^{mn}_k$ is defined by
$$
t_{k}^{mn}= (m+1)k-n, \quad 0 \leq k\leq n.
$$
 As will be seen later, for the  purpose of computing the polarization tensor, we consider only the case where $1\leq n\leq m$. 
 (If $n >m$, $u|_{\Omega}(z)$ turns out to be more complicated polynomial than $z^n$ but is still a polynomial of degree $n$.)

 Let us assume $1\leq n \leq m$. In view of \eqref{eqn_u_Omega_zn}, we define
 \begin{align*}
\ds w(\rho,\theta):&=
\begin{cases}
\ds e^{n\rho}\cos n \theta+\sum_{k=1}^n a^k\begin{pmatrix}n \\ k\end{pmatrix}e^{-t^{mn}_k\rho}\cos t^{mn}_k\theta, &\quad \rho<\rho_0,
 \\[1.5em]
\ds {e^{-n(\rho-2\rho_0)}}\cos n \theta+\sum_{k=1}^n a^k\begin{pmatrix}n \\ k\end{pmatrix}e^{-t^{mn}_k\rho}\cos t^{mn}_k\theta, &\quad \rho>\rho_0.
\end{cases}
 \end{align*} 
Note that $w$ is harmonic in $\{ \rho < \rho_0\}$ and $\{ \rho > \rho_0\}$ and $w=O(e^{-\rho})$ as $\rho\rightarrow\infty$. Moreover, 
 \begin{align}
\begin{cases}
  \ds w|_+=w|_- &\quad \mbox{ on }\rho=\rho_0,
\\[0.5em]
 \ds\frac{\p w}{\p\rho}\Big|_+ - \frac{\p w}{\p\rho}\Big|_-= 
 (-2) n e^{n\rho_0} \cos n\theta &\quad \mbox{ on }\rho=\rho_0.
 \end{cases}
\end{align}
Therefore, the function $w$ is equal to $\widetilde{u}$ up to a multiplicative constant. More precisely, we have
\be\label{eqn_tilde_u_w}
\widetilde{u}(\rho,\theta) = - \frac{1}{2n}e^{-n\rho_0} w(\rho,\theta).
\ee

Now we are ready to compute $\mathcal{K}_\Omega^*[J^{-1}\cos n\theta]$.
We can check that
\begin{align}
\frac{1}{2}\bigg( \frac{\p w}{\p\rho}\Big|^+_{\rho=\rho_0} + \frac{\p w}{\p\rho}\Big|^-_{\rho=\rho_0}\bigg)&= \sum_{k=1}^n - t^{mn}_k a^k\begin{pmatrix}n \\ k\end{pmatrix}e^{-t^{mn}_k\rho_0}\cos t^{mn}_k\theta.
 \end{align}
Then it follows from \eqref{eqn_NP_single_identity} and \eqref{eqn_tilde_u_w} that
\be\label{eqn_NP_Jcos}
\mathcal{K}_\Omega^*[J^{-1}\cos n\theta] = \frac{1}{J}\sum_{k=1}^n \delta^k\frac{ t^{mn}_k}{2n}\begin{pmatrix}n \\ k\end{pmatrix}\cos t^{mn}_k\theta
\ee
for $1\leq n\leq m$. In exactly the same manner, we can show that
\be\label{eqn_NP_Jsin}
\mathcal{K}_\Omega^*[J^{-1}\sin n\theta] = 
-\frac{1}{J}\sum_{k=1}^n \delta^k\frac{t^{mn}_k}{2n}\begin{pmatrix}n \\ k\end{pmatrix}\sin t^{mn}_k\theta.
\ee

It is worth mentioning that we can also compute the single layer potentials for $J^{-1}\cos n\theta$ and $J^{-1}\sin n\theta$:
\begin{align}
\mathcal{S}_\Omega[J^{-1}\cos n\theta] &= 
-\frac{1}{2n}\cos n \theta - \frac{1}{2n}\sum_{k=1}^n \delta^k\begin{pmatrix}n \\ k\end{pmatrix} \cos t^{mn}_k\theta, \label{SL_cos}
\end{align}
and
\begin{align}
\mathcal{S}_\Omega[J^{-1}\sin n\theta] &= 
- \frac{1}{2n}\sin n \theta+\frac{1}{2n}\sum_{k=1}^n \delta^k\begin{pmatrix}n \\ k\end{pmatrix}\sin t^{mn}_k\theta. \label{SL_sin}
\end{align}

\subsection{Asymptotic behavior of the Neumann-Poincar\'e operator $\mathcal{K}^*_\Omega$}

If $\delta$ is small enough, then the shape of $\partial \Omega$ is close to a circle.
Next we investigate the asymptotic behavior of the Neumann-Poincar\'e operator and its spectrum for small $\delta$.
From \eqref{eqn_NP_Jcos}, we infer
\begin{align}
\mathcal{K}_\Omega^*[J^{-1}\cos n\theta] 
&=\delta \frac{(m+1-n)}{2} J^{-1} \cos(m+1-n)\theta + O(\delta^2),
\nonumber
\\
\mathcal{K}_\Omega^*[J^{-1}\sin n\theta] 
&= - \delta\frac{(m+1-n)}{2} J^{-1} \sin(m+1-n)\theta + O(\delta^2)
\label{NP_order_1_m}
\end{align}
for small $\delta$ and $1\leq n\leq m$. One can verify the decay
\be\label{NP_high_order}
\mathcal{K}_\Omega^*[J^{-1}\cos n\theta],\,\mathcal{K}_\Omega^*[J^{-1}\sin n\theta] = O(\delta^2)
\ee
for small $\delta$ and $n \geq m+1$.

Let us denote by
$$v_n^c = J^{-1} \cos n\theta, \quad v_n^s = J^{-1} \sin n\theta, $$
and let $V_m^c$ and $V_m^s$ be the subspaces defined by
\begin{align*}
V_m^c = \mbox{span} \{ v_1^c,v_2^c,..., v_n^c ,..., v_m^c\}
\quad \mbox{and } 
V_m^s = \mbox{span} \{ v_1^s,v_2^s,..., v_n^s ,..., v_m^s\}.
\end{align*}
In order to better illustrate the structure of the Neumann-Poincar\'e operator, we first consider the low degree case  $m=3$.
Using $\{v_n^c\}_{n=1}^3$ as a basis, we have the following matrix representation for  $\mathcal{K}_\Omega^*$ on the subspace $V_3^c$:
\begin{align}
\mathcal{K}_\Omega^* = \frac{\delta}{2}\begin{bmatrix}
 0 & 0 & 1 \\
 0 & 2 & 0 \\
 3 & 0 & 0
 \end{bmatrix}
+O(\delta^2).
\end{align}
Similarly, using $\{v_n^s\}_{n=1}^3$ as a basis, we have the following matrix representation for  $\mathcal{K}_\Omega^*$ on the subspace $V_3^s$:
\begin{align}
\mathcal{K}_\Omega^* = - \frac{\delta}{2}\begin{bmatrix}
 0 & 0 & 1 \\
 0 & 2 & 0 \\
 3 & 0 & 0
 \end{bmatrix}
+O(\delta^2).
\end{align}

Next,  we turn to the general case. For arbitrary integer $m$, the Neumann-Poincar\'e operator $\mathcal{K}_\Omega^*$ has the following matrix representation on $V_m^c$:
\begin{align}
\mathcal{K}_\Omega^*&=\frac{\delta}{2}
\begin{bmatrix}
0 & 0 & ... &  & 1 \\
0 & ... &  & 2 &  \\
... &  &... &  & ... \\
 & m-1 &  & ... & 0 \\
m &  &... & 0 & 0 \\
\end{bmatrix} +O(\delta^2)
\nonumber
\\
:&= \frac{\delta}{2}M_{\Omega,m} + O(\delta^2).
\label{eqn_NP_mat_Vc}
\end{align}
Similarly, on $V_m^s$, we have
\begin{align}
\mathcal{K}_\Omega^*&= - \frac{\delta}{2} 
\begin{bmatrix}
0 & 0 & ... &  & 1 \\
0 & ... &  & 2 &  \\
... &  &... &  & ... \\
 & m-1 &  & ... & 0 \\
m &  &... & 0 & 0 \\
\end{bmatrix} +O(\delta^2)
\nonumber
\\
:&= - \frac{\delta}{2} M_{\Omega,m} + O(\delta^2).
\label{eqn_NP_mat_Vs}
\end{align}

Let us now consider the eigenvalues and the associated eigenvectors of the matrix $M_{\Omega,m}$. The following lemma can be easily proven.
\begin{lem}\label{lem_M_Omega_spectrum}
\begin{itemize}
\item[(i)] If $m$ is odd, that is, $m=2k-1$ for some $k\in\mathbb{N}$, then the matrix $M_{\Omega,m}$ has the following eigenvalues:
$$
k, \pm \sqrt{1 \cdot m}, \pm \sqrt{2\cdot (m-1)},..., \pm \sqrt{(k-1)\cdot (k+1)},
$$
and the associated eigenvectors are given by
$$
\mathbf{e}_{k},\quad \mathbf{e}_1\pm \sqrt{m}\,\mathbf{e}_m, \quad\mathbf{e}_2\pm\sqrt{m-1}\,\mathbf{e}_{m-1},\quad...,\quad
\sqrt{k-1}\,\mathbf{e}_{k-1}\pm\sqrt{k+1}\,\mathbf{e}_{k+1},
$$
where $\mathbf{e}_{i}$ is the unit vector in the $i$-th direction.
\item[(ii)] If $m$ is even, that is, $m=2k$ for some $k\in\mathbb{N}$,
then the matrix $M_{\Omega,m}$ has the following eigenvalues:
$$
\pm \sqrt{1 \cdot m}, \pm \sqrt{2\cdot (m-1)},..., \pm \sqrt{k\cdot (k+1)},
$$
and the associated eigenvectors are given by
$$
\mathbf{e}_1\pm\sqrt{m}\,\mathbf{e}_m, \quad \sqrt{2}\,\mathbf{e}_2\pm\sqrt{m-1}\,\mathbf{e}_{m-1},\quad...,\quad
\sqrt{k}\,\mathbf{e}_{k}\pm\sqrt{k+1}\,\mathbf{e}_{k+1}.
$$
\end{itemize}

\end{lem}

\smallskip
Using \eqref{eqn_NP_mat_Vc}, Lemma \ref{lem_M_Omega_spectrum} and the perturbation theory, we get the following asymptotic result for $\mathcal{K}^*_\Omega$ on $V_c = \mbox{span}\{ v_1^c, v_2^c, \ldots, v_m^c, \ldots \}$.
\begin{thm} \label{thm_NP_spectral_asymp_Vc}
For small $\delta$, we have the following asymptotic expansions of eigenvalues and eigenfunctions of $\mathcal{K}_\Omega^*$ on $V_c$:
\begin{itemize}
\item[(i)] If $m$ is odd, that is, $m=2k-1$ for some $k\in\mathbb{N}$:

Eigenvalues: up to order $\delta$
$$
\frac{\delta}{2}\quad\times \quad \Big\{\, k,\, \pm \sqrt{1 \cdot m},\, \pm \sqrt{2\cdot (m-1)}\,,...,\, \pm \sqrt{(k-1)\cdot (k+1)} \,\Big\}.
$$
Eigenfunctions: up to order $\delta^0$
$$
v_{k}^c,\quad v_1^c\pm \sqrt{m}\, v_m^c, \quad \sqrt{2}\,v_2^c\pm \sqrt{m-1}\,v_{m-1}^c,\quad...,\quad
\sqrt{k-1}\,v_{k-1}^c\pm \sqrt{k+1}\,v_{k+1}^c.
$$
\item[(ii)] If $m$ is even, that is, $m=2k$ for some $k\in\mathbb{N}$:

Eigenvalues: up to order $\delta$
$$
\frac{\delta}{2}\quad\times \quad \Big\{\, \pm \sqrt{1 \cdot m},\, \pm \sqrt{2\cdot (m-1)}\,,..., \,\pm \sqrt{k\cdot (k+1)} \,\Big\}.
$$
Eigenfunctions: up to order $\delta^0$
$$
v^c_1\pm \sqrt{m}\,v^c_m, \quad \sqrt{2}\,v^c_2\pm \sqrt{m-1}\,v^c_{m-1},\quad...,\quad
\sqrt{k}\,v_k^c\pm \sqrt{k+1}\,v^c_{k+1}.
$$
\end{itemize}

\end{thm}

Similarly, we have the following result for $\mathcal{K}_\Omega^*$ on the subspace $V_s = 
\mbox{span}\{ v_1^s, v_2^s, \ldots, v_m^s, \ldots \}$.

\begin{thm}  \label{thm2.3}
We have the following asymptotic expansion of eigenvalues and eigenfunctions of the Neumann-Poincar\'e operator $\mathcal{K}_\Omega^*$ on the subspace $V_s$ for small $\delta$:
\begin{itemize}
\item[(i)] If $m$ is odd, that is, $m=2k-1$ for some $k\in\mathbb{N}$:

Eigenvalues: up to order $\delta$
$$
- \frac{\delta}{2}\quad\times \quad \bigg\{
  k,\, \pm \sqrt{1 \cdot m},\, \pm \sqrt{2\cdot (m-1)}\,,...,\, \pm \sqrt{(k-1)\cdot (k+1)} \bigg\}.
$$
Eigenfunctions: up to order $\delta^0$
$$
v_{k}^s,\quad v_1^s\pm \sqrt{m}v_m^s, \quad \sqrt{2}v_2^s\pm \sqrt{m-1}v_{m-1}^s,\quad...,\quad
\sqrt{k-1}\,v_{k-1}^s\pm \sqrt{k+1}\,v_{k+1}^s.
$$
\item[(ii)] If $m$ is even, that is, $m=2k$ for some $k\in\mathbb{N}$:

Eigenvalues: up to order $\delta$
$$
- \frac{\delta}{2}\quad\times \quad \bigg\{
\pm \sqrt{1 \cdot m},\, \pm \sqrt{2\cdot (m-1)}\,,..., \,\pm \sqrt{k\cdot (k+1)} \bigg\}.
$$
Eigenfunctions: up to order $\delta^0$
$$
v^s_1\pm \sqrt{m}\,v^s_m, \quad \sqrt{2}\,v^s_2\pm \sqrt{m-1}\,v^s_{m-1},\quad...,\quad
\sqrt{k}\,v^s_k\pm \sqrt{k+1}\,v^s_{k+1}.
$$
\end{itemize}

\end{thm}

\begin{cor} \label{cor2.4}
Suppose that $m$ is odd, that is, $m=2k-1$ for some $k\in\mathbb{N}$. In other words, 
$\Omega$ is a star-shaped domain with $2k$ petals. Then, up to order $\delta$, the Neumann-Poincar\'e operator $\mathcal{K}_\Omega^*$ has the following $2k$ eigenvalues:
$$
\frac{\delta}{2}\quad\times \quad \bigg\{ \pm \sqrt{1 \cdot m},\, \pm \sqrt{2\cdot (m-1)}\,,...,\, \pm \sqrt{(k-1)\cdot (k+1)}, \,\pm\sqrt{k\cdot k} \bigg\}.
$$
\end{cor}

\subsection{Generalized polarization tensors and their spectral representations}

\subsubsection{First-order polarization tensor}
Let us compute the first-order polarization tensor associated with  $\Omega$ and $\lambda$. 
For simplicity, we consider only the case when $m$ is odd, that is, $m=2k-1$ for some $k\in\mathbb{N}$. 
The case where $m$ is even can be treated analogously. Numerical results are presented in section \ref{sect4} for both cases.  

Since $m$ is odd, the shape of $\Omega$ has even symmetry with respect to both $x_1$-axis and $x_2$-axis. Thanks to this symmetry,  $M(\lambda,\Omega)$ has the following simple form \cite{book2}:
$$
M(\lambda,\Omega) = m_{11} \begin{bmatrix}1&0\\ 0 &1 \end{bmatrix},
$$
where $m_{11}$ is given by
\begin{align*}
m_{11} = \left\langle x_1, (\lambda I-\mathcal{K}^*_\Omega)^{-1} [\nu_1] \right\ra_{L^2}.
\end{align*}
Let $\lambda_n$ and $\varphi_n$, $n \in \mathbb{N}$,  be the eigenvalues and the (normalized) eigenfunctions of $\mathcal{K}^*_\Omega$, respectively.
Then, from the spectral decomposition of $\mathcal{K}^*_\Omega$, we have \cite{matias, kang1, hyeonbae}
\begin{align*}
m_{11}&=  \sum_{j} \frac{1}{\lambda-\lambda_{j}}\frac{\la x_1,\varphi_{j}\ra_{L^2}\la \varphi_{j}, -\mathcal{S}_\Omega[\nu_1] \ra_{L^2}}{ \la \varphi_j , -\mathcal{S}_\Omega[\varphi_j] \ra_{L^2} }
\\
&=\sum_{j} \frac{(\frac{1}{2}-\lambda_{j})}{\lambda-\lambda_{j}}\frac{|\la x_1,\varphi_{j}\ra_{L^2(\p\Omega)}|^2
}{\la \varphi_j , - \mathcal{S}_\Omega[\varphi_j] \ra_{L^2}}.
\end{align*}

By Theorems \ref{thm_NP_spectral_asymp_Vc} and \ref{thm2.3},  one can see that only the following two eigenvalues and two eigenfunctions contribute to $m_{11}$ up to order $\delta$:
\begin{align*}
&\mbox{eigenvalues } \lambda_\pm:=\pm  \frac{1}{2}\delta \sqrt{m}, 
\\ 
&\mbox{eigenfunctions } \varphi_\pm:=v_1^c\pm \sqrt{m} \,v_m^c .
\end{align*}
In fact, for other eigenfunctions, we have $\langle x_1, \varphi_j\rangle_{L^2} = O(\delta)$. In what follows we calculate $\langle x_1, \varphi_\pm\rangle_{L^2}$ and $\langle x_1, \mathcal{S}_\Omega[\varphi_\pm]\rangle_{L^2}$.  

First, since $d\sigma = J d\theta$ and 
$$
x_1|_{\p \Omega} = \Re\{ \Phi(e^{\rho_0+i\theta})\} = e^{\rho_0}\cos\theta + a e^{-m\rho_0}\cos m\theta,
$$
we have
\begin{align*}
\la x_1,\varphi_\pm\ra_{L^2(\p\Omega)} &= \int_{0}^{2\pi} (e^{\rho_0}\cos\theta + a e^{-m\rho_0}\cos m\theta) (\cos\theta\pm \sqrt{m} \cos m\theta) \,d\theta, 
\\
&=\pi e^{\rho_0}(1 + 2 \lambda_\pm).
\end{align*}
Now, we compute $\la \varphi_\pm, \mathcal{S}_\Omega[\varphi_\pm]\ra_{L^2}$. Note that, from \eqnref{SL_cos}, we have
$$
\mathcal{S}_\Omega[v_n^c] = 
- \frac{1}{2n}\cos n \theta
- \frac{1}{2}\delta \cos (m+1-n)\theta 
+O(\delta^2).
$$
Consequently,
\begin{align*}
\la \varphi_\pm, -\mathcal{S}_\Omega[\varphi_\pm] \ra_{L^2} &= 
\la v_1^c \pm \sqrt{m}\,v_m^c, -\mathcal{S}_{\Omega}[v_1^c \pm \sqrt{m}\,v_m^c]\ra_{L^2(\p\Omega)}
\\
&=
\int_{0}^{2\pi} (\cos\theta \pm  \sqrt{m}\,\cos m\theta) 
\\
&\qquad\times
(\frac{1}{2}\cos\theta 
+\frac{\delta}{2} \cos m\theta
\pm\frac{1}{2\sqrt m}\cos m\theta 
\pm\sqrt m\,\frac{\delta}{2} \cos \theta
) \,d\theta
+O(\delta^2)
\\
&=\pi(1 +2\lambda_\pm) +O(\delta^2).
\end{align*}

Finally,  we are ready to obtain an approximation formula for $m_{11}$. 

\begin{thm}
We have
\beq\label{PT_approx_m11}
m_{11} =
 \frac{\pi}{2} e^{2\rho_0}
 \Big(\frac{1}{\lambda-\lambda_{+}}  +  \frac{1}{\lambda-\lambda_{-}}\Big) + O(\delta^2),
\eeq
as $\delta \rightarrow 0$. 
\end{thm}

\subsubsection{Second-order contracted generalized polarization tensors}

Let $M_{mn}^{cc},M_{mn}^{ss},M_{mn}^{sc},$ and $M_{mn}^{cs}$ be the contracted generalized polarization tensors.
One can easily see that $M_{22}^{sc}=M_{22}^{cs}=0$ and $M_{12}=M_{21}=0$. 
 We only need to consider $M_{22}^{cc}$ and $M_{22}^{ss}$. 
It turns out that only the following two eigenvalues and two eigenfunctions contribute to $M_{22}^{cc}$(up to the order $\delta$):
\begin{align*}
&\mbox{eigenvalues } \lambda_\pm':=\pm  \frac{1}{2}\delta \sqrt{2\cdot(m-1)},
\\ 
&\mbox{eigenfunctions } \varphi_\pm':=\sqrt{2}\,v_2^c\pm \sqrt{m-1} \,v_{m-1}^c.
\end{align*}
Let $H:=\Re\big\{ (x_1+i x_2)^2 \big\}$. Then we have  
$$H|_{\p\Omega}= e^{2\rho_0}(\cos 2\theta+2\delta \cos m\theta+\delta^2\cos 2m\theta).$$
Therefore, 
\begin{align*}
\la H,\varphi_\pm'\ra_{L^2(\p\Omega)} &= \int_{0}^{2\pi} e^{2\rho_0}(\cos 2\theta+2\delta \cos m\theta+\delta^2\cos 2m\theta) (\sqrt{2}\,\cos 2\theta\pm \sqrt{m-1} \cos (m-1)\theta) \,d\theta
\\
&=\sqrt{2}\pi e^{2\rho_0}.
\end{align*}
Now we compute $\la \varphi'_\pm, - \mathcal{S}_\Omega[\varphi'_\pm]\ra_{L^2}$. Since 
$$
\mathcal{S}_\Omega[v_n^c] = 
-\frac{1}{2n}\cos n \theta
-\frac{1}{2}\delta \cos (m+1-n)\theta 
+O(\delta^2),
$$
we obtain
\begin{align*}
\la \varphi_\pm', - \mathcal{S}_\Omega[\varphi_\pm']\ra_{L^2} &= 
\la \sqrt{2}\,v_2^c \pm \sqrt{m-1}\,v_{m-1}^c, - \mathcal{S}_{\Omega}[\sqrt{2}\,v_2^c \pm \sqrt{m-1}\,v_{m-1}^c]\ra_{L^2(\p\Omega)}
\\
&= 
\int_{0}^{2\pi} (\sqrt{2}\,\cos 2\theta \pm  \sqrt{m-1}\,\cos (m-1)\theta) 
\\
&\quad\times
(\frac{1}{2\sqrt{2}}\cos 2\theta 
+\frac{\delta}{\sqrt{2}} \cos (m-1)\theta\\ 
&\quad
\pm\frac{1}{2\sqrt {m-1}}\cos (m-1)\theta 
\pm\sqrt {m-1}\,\frac{\delta}{2} \cos 2\theta
) \,d\theta
+O(\delta^2)
\\
&=  \pi(1 +2\lambda_\pm') +O(\delta^2).
\end{align*}
Finally, we find
\begin{align}
M_{22}^{cc}
&= \sum_{j} \frac{(\frac{1}{2}-\lambda_{j})}{(\lambda-\lambda_{j})}\frac{|\la H,\varphi_{j}\ra_{L^2(\p\Omega)}|^2
}{\la \varphi_j , - \mathcal{S}_\Omega[\varphi_j] \ra_{L^2}},
\nonumber
\\
&=
{\pi} e^{4\rho_0} \Big(
\frac{(\frac{1}{2}-\lambda_+')}{(\frac{1}{2}+\lambda_+')
(\lambda-\lambda_+')} + \frac{(\frac{1}{2}+ \lambda_-')}{(\frac{1}{2} -\lambda_-')(\lambda-\lambda_-')} 
\Big)
+O(\delta^2).\label{CGPT_approx_22}
\end{align}
Similarly, one can show that $M_{22}^{ss}$ has similar asymptotic expansion.

\subsection{Classification of algebraic domains in the class $\mathcal{Q}$}
The identification of the parameters $\rho_0$ and $m$ is now straightforward using the results of the previous subsection. 
Suppose that we can obtain the values of $\lambda_\pm, \lambda_\pm'$ approximately from $m_{11}$ and $M_{22}^{cc}$. Then, by formula \eqnref{PT_approx_m11}, we can easily find the parameter $\rho_0$, which determines the size of $\Omega$. In order to reconstruct the parameters $m$ and $\delta$ we turn to the following equations:
$$
\lambda_+ = \frac{1}{2}\delta \sqrt{m}, \quad \lambda_+' = \frac{1}{2}\delta \sqrt{2\cdot (m-1)}.
$$
Solving the above equations for $m$ and $\delta$ yields:
$$
 m = \frac{ \lambda_+^2 }{  \lambda_+^2 -(\lambda_+')^2/2 }
 ,\quad
 \delta = {2}  {  \sqrt{  \lambda_+^2 - (\lambda_+')^2/2  } }.
$$

\section{Plasmonic resonances for two separated disks} \label{sect3}
In this section, we consider the spectrum of the Neumann-Poincar\'e operator when two conductors are  located closely to each other in $\RR^2$. As an application of the spectral decomposition of the Neumann-Poincar\'e operator, we derive  the $(1,1)$-entry, $m_{11}$, of the first-order polarization tensor associated with the two disks.

\subsection{The bipolar coordinates and the boundary integral operators}

Let $B_1$ and $B_2$ be two disks with conductivity $k$ embedded in the background with conductivity $1$. The conductivity $k$ is such that  $0<k\neq1<\infty$. Let $\sigma_{B_1 \cup B_2}$ denote the conductivity distribution, \textit{i.e.},
\beq\label{sigma:def}
\sigma_{B_1 \cup B_2}=k\chi(B_1)+k\chi(B_2)+\chi(\RR^2\setminus (B_1\cup B_2),
\eeq
where $\chi$ is the characteristic function.
Let $\ep$ be the distance between two disks, that is,
$$\ep:=\mbox{dist}(B_1,B_2).$$
 We set Cartesian coordinates $(x_1,x_2)$ such that $x_1$-axis is parallel to the line joining the centers of the two disks.\\

\textbf{(Definition)} Each point $\mathbf{x}=(x,y)$ in the Cartesian coordinate system corresponds to  $(\xi,\theta)\in\RR\times (-\pi,\pi]$ in the bipolar coordinate system through the equations
\begin{equation} \label{bipolar}
x=\alpha\frac{ \sinh \xi }{\cosh \xi - \cos \theta} \quad \mbox{ and } \quad y=\alpha\frac{\sin \theta}{\cosh \xi - \cos \theta}
\end{equation}
with a positive number $\alpha$. In fact, the bipolar coordinates can be defined using a conformal mapping. Define a conformal map $\Psi$ by
$$
z=x+i y= \Psi(\zeta)=\alpha \frac{\zeta+1}{\zeta-1}.
$$
If we write $\zeta=e^{\xi- i\theta}$, then we can recover \eqnref{bipolar}.\\

\noindent\textbf{(The coordinate curve)} From the definition, we can derive that the coordinate curves  $\{\xi=c\}$ and $\{\theta=c\}$ are, respectively, the zero-level set of the following two functions:
\beq\label{function:f}f_\xi(x,y)=\left(x-\alpha\frac{\cosh c}{\sinh c}\right)^2 +y^2-\left(\frac{\alpha}{\sinh c}\right)^2\eeq
and $$f_{\theta}(x,y)=x^2 +\left(y-\alpha\frac{\cos c}{\sin c}\right)^2-\left(\frac{\alpha}{\sin c}\right)^2.$$\\

\noindent\textbf{(Basis vectors)} Orthonormal basis vectors $\{\hat{\mathbf{e}}_{\xi},\hat{\mathbf{e}}_{\theta}\}$ are defined as follows:
$$
\hat{\mathbf{e}}_{\xi}:= \frac{\p \mathbf{x}/ \p \xi}{|\p \mathbf{x}/ \p\xi|} \quad \mbox{and} \quad \hat{\mathbf{e}}_{\theta}:= \frac{\p \mathbf{x}/ \p \theta}{|\p \mathbf{x}/ \p\theta|}.
$$\\

\noindent \textbf{(Normal- and tangential derivatives and line element)} 
In the bipolar coordinates, the scaling factor $h$ is 
$$
h(\xi,\theta):=\frac{\cosh\xi-\cos\theta}{\alpha}.
$$
The gradient of any scalar function $g$ is
\beq \label{grad_bipolar}
\nabla g = h(\xi,\theta)\left( \frac{\p g}{\p\xi}\hat{\mathbf{e}}_{\xi}+ \frac{\p g}{\p\theta}\hat{\mathbf{e}}_{\theta}\right).
\eeq
Moreover, the normal and tangential derivatives of a function $u$ in bipolar coordinates are 
\begin{eqnarray}
\left\{ \begin{array}{l} \label{nor_bipolar}
\ds \pd{u}{\nu}\Bigr|_{\xi=c}=\nabla u\cdot v_{\xi=c}=-\mbox{sgn}(c)h(c,\theta)\pd{u}{\xi}\Bigr|_{\xi=c},\\
\nm
 \label{tan_bipolar} \ds \pd{u}{T}\Bigr|_{\xi=c}=-\mbox{sgn}(c)h(c,\theta)\pd{u}{\theta}\Bigr|_{\xi=c},
\end{array} \right. \end{eqnarray}
and the line element $d\sigma$ on the boundary $\{\xi=\xi_0\}$ is
$$
d\sigma = \frac{1}{h(\xi_0,\theta)} d\theta.\\
$$

\noindent \textbf{(Separation of variables)} The bipolar coordinate system admits separation of variables for any harmonic function $f$ as follows:
\begin{align}
f(\xi,\theta)&=a_0+b_0\xi+c_0\theta+\sum_{n=1}^\infty\bigr[ (a_n e^{n\xi}+b_n e^{-n\xi})\cos n\theta+\bigr (c_n e^{n\xi}+d_n e^{-n\xi})\sin n\theta\bigr],
\end{align}
where $a_n$, $b_n$, $c_n$ and $d_n$ are constants.

For $\xi>0$, we have
\begin{align} \label{xy_zeta}
\frac{\sinh \xi -i\sin \theta}{\cosh \xi -\cos\theta} &=\frac{e^\zeta +e^{-\zeta}}{e^\zeta - e^{-\zeta}}=1+2\sum_{n=1}^\infty e^{-n\xi}(\cos n\theta-i\sin n\theta),
\end{align}
with $\zeta=({\xi+i\theta})/{2}$.

Using \eqnref{bipolar}, we have the following harmonic expansions for the two linear functions $x_1$ and $x_2$:
\beq\label{x1_bipolar}
x_1=\mbox{sgn}({\xi})\alpha\left[1+2\sum_{n=1}^\infty  e^{-n|\xi|}\cos n\theta\right],
\eeq
and
$$x_2=2\alpha\sum_{n=1}^\infty  e^{-n|\xi|}\sin n\theta.$$

Let $\mathbb{K}^*$ be the Neumann-Poincar\'{e} operator given by
$$\mathbb{K}^*:=
\left[
  \begin{array}{cc}
  \ds  \Kcal^*_{B_1} & \ds \pd{}{\nu^{(1)}}\Scal_{B_2} \\
  \ds \pd{}{\nu^{(2)}}\Scal_{B_1} & \ds\Kcal_{B_2}^* \\
   \end{array}
\right],$$
and define the operator $\mathbb{S}$ by
$$\mathbb{S}=
\left[\begin{array}{cc}
\Scal_{B_1}&\Scal_{B_2}\\
\Scal_{B_1}&\Scal_{B_2}
\end{array}\right].$$
Here, $\nu^{(i)}$ is the outward normal on $\partial B_i$, $i=1,2$. 
 
Then, from \cite{arma}, $\mathbb{K}^*$ is self-adjoint with the inner product
\beq
\la \varphi,\psi\ra_{\mathcal{H}}:=-\la\varphi,\mathbb{S}[\psi]\ra,\quad\mbox{for }\varphi,\psi\in L^2(\partial B_1)\times L^2(\partial B_2).\eeq
\subsection{Neumann Poincar\'e-operator for two separated disks and its spectral decomposition}\label{NP_separated}
First we introduce some notations. Set
\beq\label{def_alpha_xi0}
\alpha= \sqrt{ \ep (r+ \frac{\ep}{4})}\quad\mbox{and}
\quad\xi_0=\sinh^{-1}\left(\frac{\alpha}{r}\right),\quad \mbox{for }j=1,2,\eeq
where $r$ is the radius of the two disks and $\ep$ their separation distance. Note that
\beq
\p B_j=\{\xi=(-1)^j\xi_0\},\quad \mbox{for }j=1,2.
\eeq

Let us denote the Neumann-Poincar\'e operator for two disks separated by a distance $\ep$ by $\mathbb{K}_\ep^*$. To find out the spectral decomposition of the Neumann-Poincar\'e operator $\mathbb{K}_\ep^*$, we use the following lemma \cite{arma}.
\begin{lem}\label{trans_eigen_lemma}
Assume that there exists  $u$ a nontrivial solution to the following equation:
\beq\label{transmission_eigen}
\begin{cases}
\ds\Delta u = 0 \quad& \mbox{in } B_1 \cup B_2\cup \mathbb{R}^2\setminus \overline{(B_1\cup B_2)},\\
\ds u|_+=u|_- \quad& \mbox{on } \p B_j,  j=1,2,\\
\ds \frac{\p u}{\p \nu}\Big|_+=k\frac{\p u}{\p \nu}\Big|_- \quad&\mbox{on } \p B_j, j=1,2,\\
\ds u(\Bx)\rightarrow 0 \quad& \mbox{as } |\Bx| \rightarrow \infty,
\end{cases}
\eeq
where $\ds k=-\frac{1+2\lambda}{1-2\lambda}<0$.
If we set $$\psi_j:=\frac{\p u}{\p \nu}\Big|^+_{\p B_j}-\frac{\p u}{\p \nu}\Big|_{\p B_j}^-, \quad \mbox{for }j=1,2,$$
then $\psi=\left[\begin{array}{c} \psi_1 \\ \psi_2 \end{array} \right]$ is an eigenvector of $\mathbb{K}^*_\ep$ corresponding to the eigenvalue $\lambda$.  
\end{lem}
One can see that the following function $u_n$ is a solution to \eqnref{transmission_eigen}:
\beq\label{un_pm}
u_n^\pm(\xi,\theta)=(\mbox{const.})+
\begin{cases}
\ds \mp \frac{1}{2|n|}(e^{|n|\xi_0}\mp e^{-|n|\xi_0}) e^{|n|\xi+in\theta},\quad&\mbox{for }\xi<-\xi_0,\\[3mm]
\ds \frac{1}{2|n|}e^{-|n|\xi_0}(e^{|n|\xi}\mp e^{-|n|\xi})e^{in\theta},\quad&\mbox{for }-\xi_0<\xi<\xi_0,\\[3mm]
\ds \frac{1}{2|n|}(e^{|n|\xi_0}\mp e^{-|n|\xi_0}) e^{-|n|\xi+in\theta},\quad&\mbox{for }\xi>\xi_0.
\end{cases}
\eeq
From \eqnref{un_pm} and Lemma \ref{trans_eigen_lemma}, we obtain eigenvalues and
eigenvectors to $\mathbb{K}_\ep^*$
$$
\lambda^\pm_{\ep,n} = \pm\frac{ 1}{2}e^{-2|n|\xi_0} \quad\mbox{and}\quad\Phi_{\ep,n}^{\pm}(\theta)= e^{i n\theta}\left[\begin{array}{c}
h(-\xi_0,\theta)\\
\mp h(\xi_0,\theta)\end{array}\right].
$$
Note that the above eigenvectors are not normalized. 

We compute $\la \Phi_{\ep,n}^\pm, -\mathbb{S}[\Phi_{\ep,n}^\pm] \ra_{L^2}$.
From \eqnref{un_pm}, one can see that
$$
\mathbb{S} [\Phi_{\ep,n}^\pm]=(\mbox{const.})+
\left[\begin{array}{c}
\mp\frac{1}{2|n|}(1\mp e^{-2|n|\xi_0})e^{in\theta}\\
\frac{1}{2|n|}(1\mp e^{-2|n|\xi_0})e^{in\theta}\end{array}\right].
$$
It follows that
$$
\la \Phi_{\ep,n}^\pm, -\mathbb{S}[\Phi_{\ep,n}^\pm] \ra_{L^2} = \frac{2\pi}{|n|}(1\mp e^{-2|n|\xi_0}).
$$
Therefore, we arrive at the following result. 
\begin{thm}
 We have the following spectral decomposition of $\mathbb{K}_\ep^*$:
\beq
\mathbb{K}_\ep^* = \sum_{n\neq 0} \frac{1}{2}e^{-2|n|\xi_0} \Psi_{\ep,n}^+ \otimes \Psi_{\ep,n}^+
+ \sum_{n\neq 0} \left( -\frac{1}{2}e^{-2|n|\xi_0}\right) \Psi_{\ep,n}^- \otimes \Psi_{\ep,n}^-,
\eeq
where $\Psi_{\ep,n}^\pm$ are the normalized eigenvectors defined by
\beq\label{normal_eigenvec}
\Psi_{\ep,n}^\pm (\theta):= \frac{\sqrt{|n|}e^{in\theta}}{\sqrt{2\pi (1\mp e^{-2|n|\xi_0})}}
\left[\begin{array}{c}
h(-\xi_0,\theta)\\
\mp h(\xi_0,\theta)\end{array}\right].
\eeq
\end{thm}
Note that
\begin{align}
(\mathcal{S}_{B_1}[\Psi_{\ep,n,1}^\pm]+\mathcal{S}_{B_2}[\Psi_{\ep,n,2}^\pm])(\xi,\theta)&=(\mbox{const.})+
 \frac{\sqrt{|n|}}{\sqrt{2\pi (1\mp e^{-2|n|\xi_0})}}
 \\
 & \quad\times
\begin{cases}
\ds \mp \frac{1}{2|n|}(e^{|n|\xi_0}\mp e^{-|n|\xi_0}) e^{|n|\xi+in\theta},\quad&\mbox{for }\xi<-\xi_0,\\[3mm]
\ds \frac{1}{2|n|}e^{-|n|\xi_0}(e^{|n|\xi}\mp e^{-|n|\xi})e^{in\theta},\quad&\mbox{for }-\xi_0<\xi<\xi_0,\\[3mm]
\ds \frac{1}{2|n|}(e^{|n|\xi_0}\mp e^{-|n|\xi_0}) e^{-|n|\xi+in\theta},\quad&\mbox{for }\xi>\xi_0.
\end{cases}
\label{Single_normal_eigen}
\end{align}

\subsection{The Polarization tensor}
Let us compute the $(1,1)$-entry $m^\ep_{11}$ of the first-order  polarization tensor for two separated disks. Our approach here is different from \cite{sinum}. It is based on the spectral decomposition of $\mathbb{K}_\ep^*$. Note that
\begin{align*}
m^\ep_{11} = \left\la \varphi, (\lambda \mathbb{I} -\mathbb{K}_\ep^*)^{-1} [\psi] \right\ra_{L^2},
\end{align*}
where
$$
\phi=\left[\begin{array}{c}
x_1|_{\p B_1}\\
x_1|_{\p B_2} \end{array}\right],\quad
\psi=\left[\begin{array}{c}
\nu_1|_{\p B_1}\\
\nu_1|_{\p B_2}\end{array}\right].
$$
The spectral decomposition of $\mathbb{K}_\ep^*$ implies
\begin{align*}
m^\ep_{11}&= \sum_{n \neq 0} \frac{\la\phi,\Psi_{\ep,n}^+\ra_{L^2}\la \Psi_{\ep,n}^+,\psi \ra_{\mathcal{H}}}{\lambda-\lambda_{\ep,n}^+}
+ \sum_{n \neq 0} \frac{\la\phi,\Psi_{\ep,n}^-\ra_{L^2}\la \Psi_{\ep,n}^-,\psi \ra_{\mathcal{H}}}{\lambda-\lambda_{\ep,n}^-}\\
&=\sum_{n \neq 0} \frac{\left(\frac{1}{2}-\lambda_{\ep,n}^+\right)|\la\phi,\Psi_{\ep,n}^+\ra_{L^2}|^2
}{\lambda-\lambda_{\ep,n}^+}
+ \sum_{n \neq 0} \frac{\left(\frac{1}{2}-\lambda_{\ep,n}^-\right)|\la\phi,\Psi_{\ep,n}^-\ra_{L^2}|^2}{\lambda-\lambda_{\ep,n}^-}.
\end{align*}
From \eqnref{x1_bipolar}, we derive the expansion
\beq\label{x1_bipolar_2}
x_1 = \mbox{sgn}(\xi)  \alpha \sum_{m=-\infty}^\infty  e^{-|m||\xi|+im\theta}.
\eeq
Therefore,
\begin{align*}
\la\phi,\Psi_{\ep,n}^+\ra_{L^2} &=2
\int_0^{2\pi} \left[-\alpha \sum_{m=-\infty}^\infty  e^{-|m|\xi_0+im\theta}\right] \frac{\sqrt{|n|}h(-\xi_0,\theta)e^{-in\theta}}{\sqrt{2\pi (1- e^{-2|n|\xi_0})}} \frac{1}{h(-\xi_0,\theta)}d\theta \\
&=-2 \sqrt{2\pi}\alpha \frac{\sqrt{|n|}e^{-|n|\xi_0} }{\sqrt{1- e^{-2|n|\xi_0}}},
\end{align*}
and
$$
\la\phi,\Psi_{\ep,n}^-\ra_{L^2}=0.
$$
As a consequence, we arrive at the following result.
\begin{prop}
We have
$$
m^\ep_{11}=\sum_{n \neq 0} \frac{4\pi \alpha^2 |n| e^{-2|n|\xi_0} }{\lambda-\lambda_{\ep,n}^+}
= 8\pi \alpha^2 \sum_{n=1}^\infty \frac{n e^{-2n\xi_0}}{\lambda-\frac{1}{2}e^{-2n\xi_0}},
$$
where $\alpha$ is given by (\ref{def_alpha_xi0}). 
\end{prop}

\subsection{Reconstruction of the separation distance}

Suppose that the first eigenvalue 
$$ 
\lambda_{\ep,1}^+ = \frac{1}{2}e^{-2\xi_0}
$$
is measured. Then we immediately find the value of $e^{\xi_0}$. From \eqnref{def_alpha_xi0}, we have
\beq
r \cosh \xi_0= \frac{\ep}{2} + r.
\eeq
By solving the above quadratic equation, we can determine 
the distance $\ep$ between the two disks.

\section{Numerical illustrations} \label{sect4}
In this section we illustrate our main findings in this paper with several numerical examples.

We use the material parameters of gold nanoparticles and suppose that we can measure their first- and second-order polarization tensors for a range of wavelengths in the visible regime. 

\begin{figure}[!h]
	\begin{center}
		\includegraphics[scale=0.17]{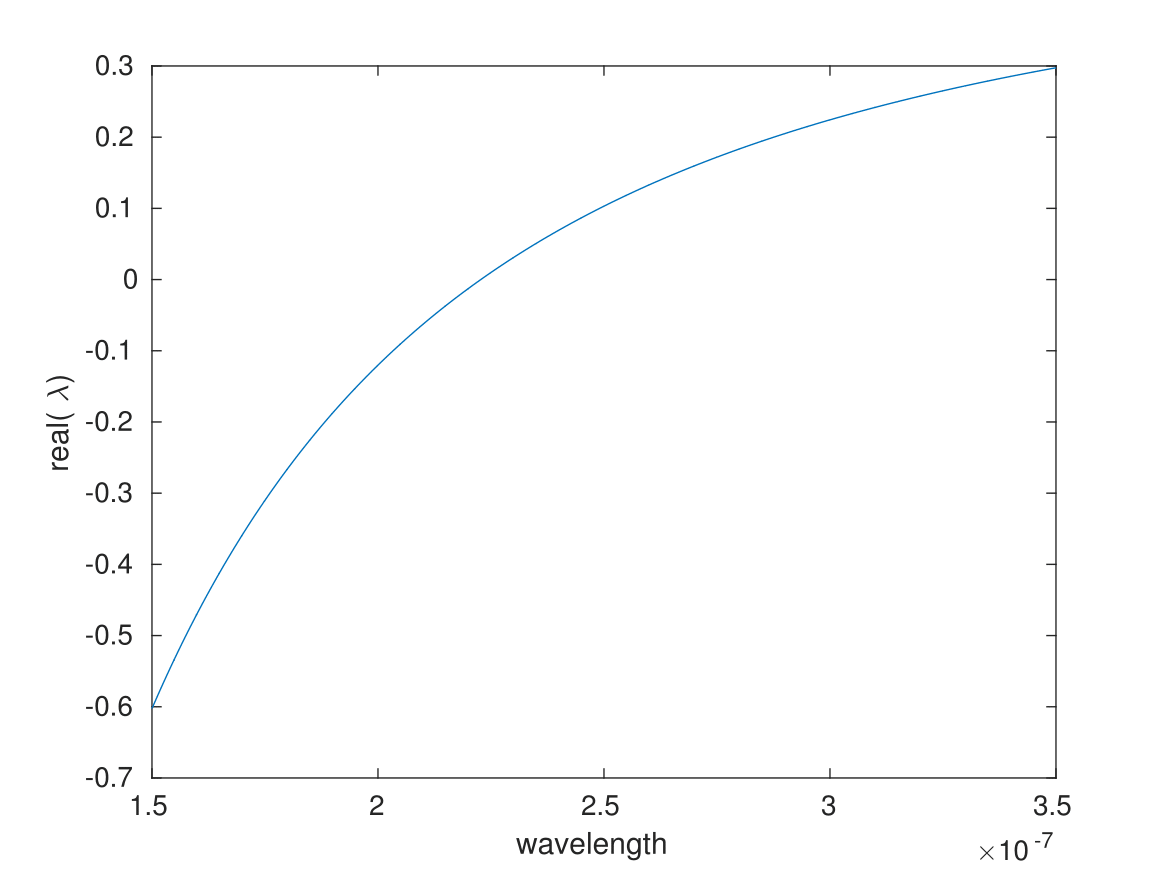}
		\includegraphics[scale=0.17]{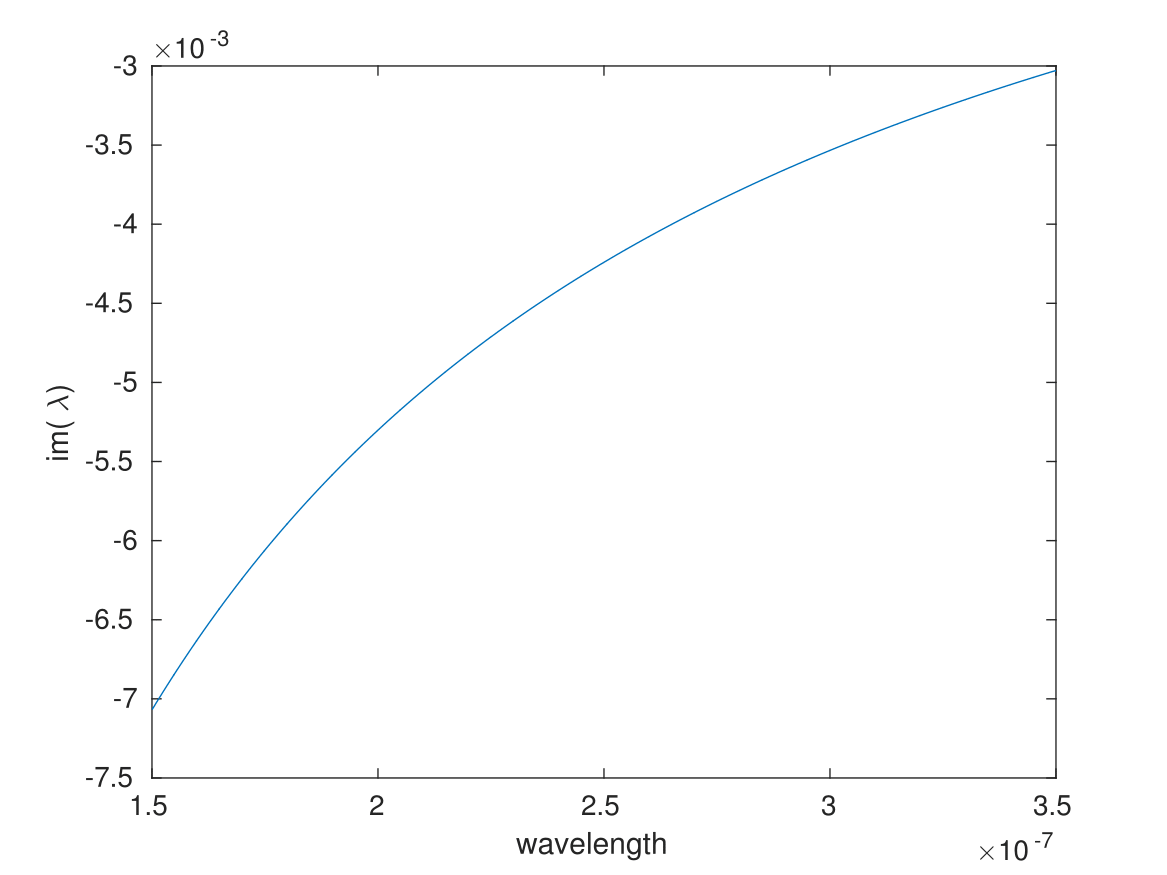}
		\caption{ \label{fig-lambda algebraic} 
		Real and imaginary parts of $\lambda$ as  function of the wavelength.}
	\end{center}
\end{figure}

Figure \ref{fig-lambda algebraic} shows the variations of the real and imaginary parts of $\lambda$, defined by (\ref{deflambda}),  as  function of the wavelength using Drude's model for $\sigma =\sigma(\omega)$, which is depending on the operating frequency $\omega$ \cite{pierre}.

As shown in Figure \ref{fig-lambda algebraic}, the imaginary part of $\lambda$ is very small. Therefore, when the real part of $\lambda$ hits an eigenvalue that contributes to the first-order polarization tensor (and therefore to the plasmonic resonances), we should see a peak in the graph of $|m_{11}|$ and $|M_{22}^{cc}|$ with respect to the wavelength. This allow us, in the case of class $\mathcal{Q}$ of algebraic domains, to recover $\lambda_+$ and $\lambda_+'$ and, in the case of two separated disks, to recover $\lambda_{\epsilon,1}^+$.

\begin{figure}[!h]
	\begin{center}
		\includegraphics[scale=0.15]{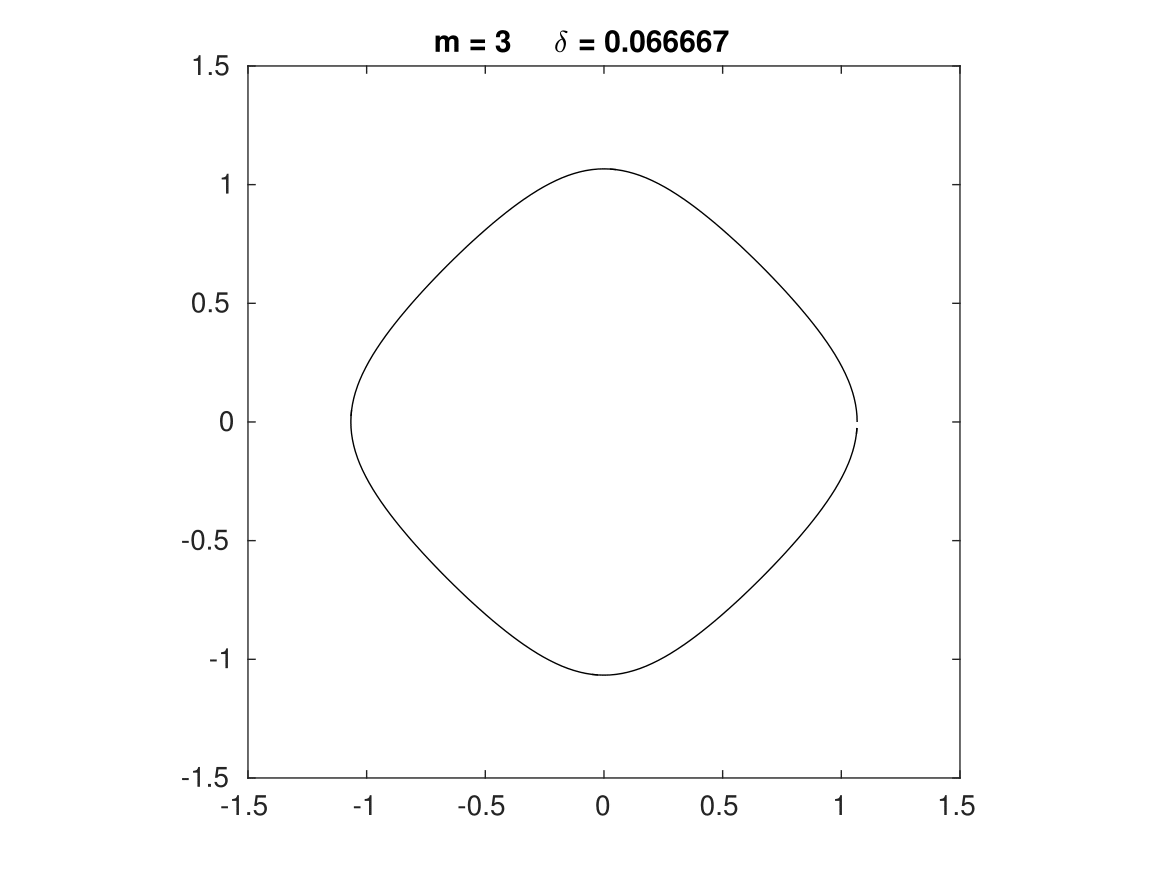}
		\includegraphics[scale=0.15]{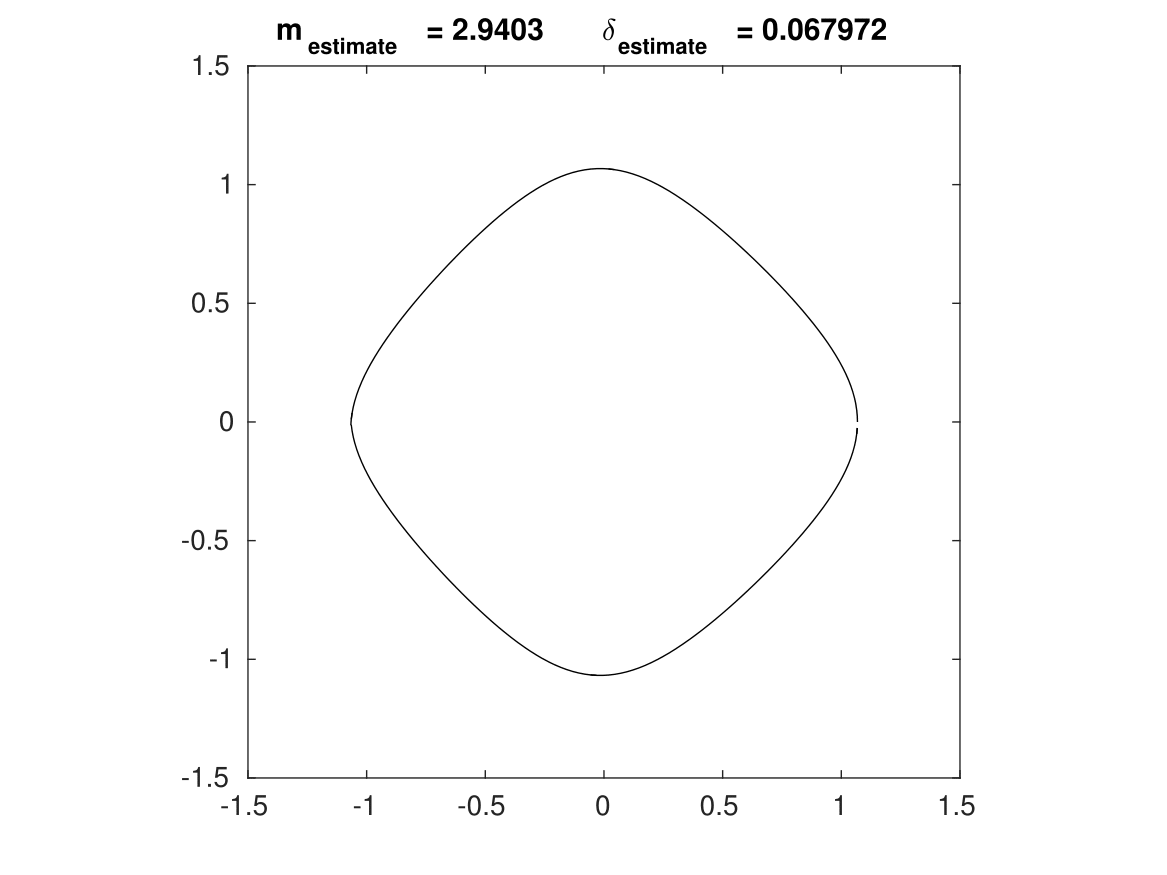}
		\includegraphics[scale=0.15]{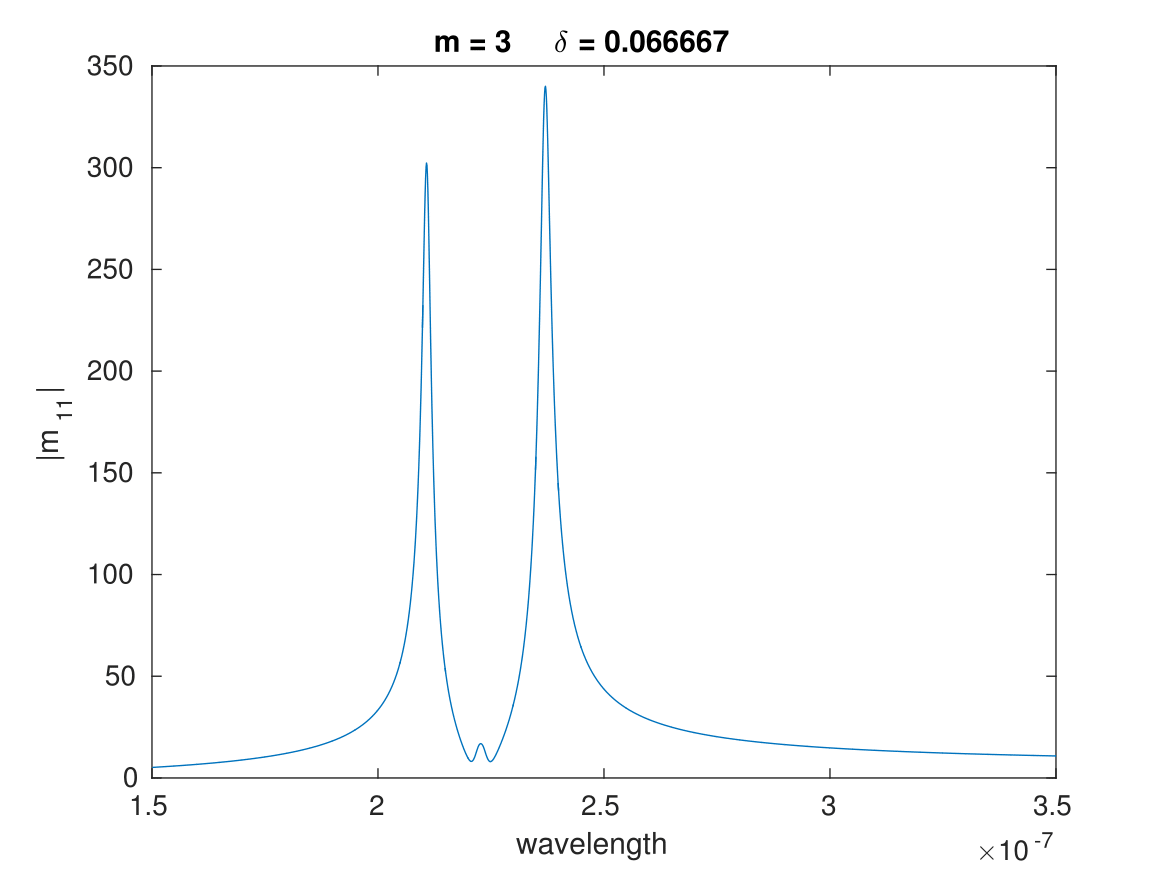}
		\includegraphics[scale=0.15]{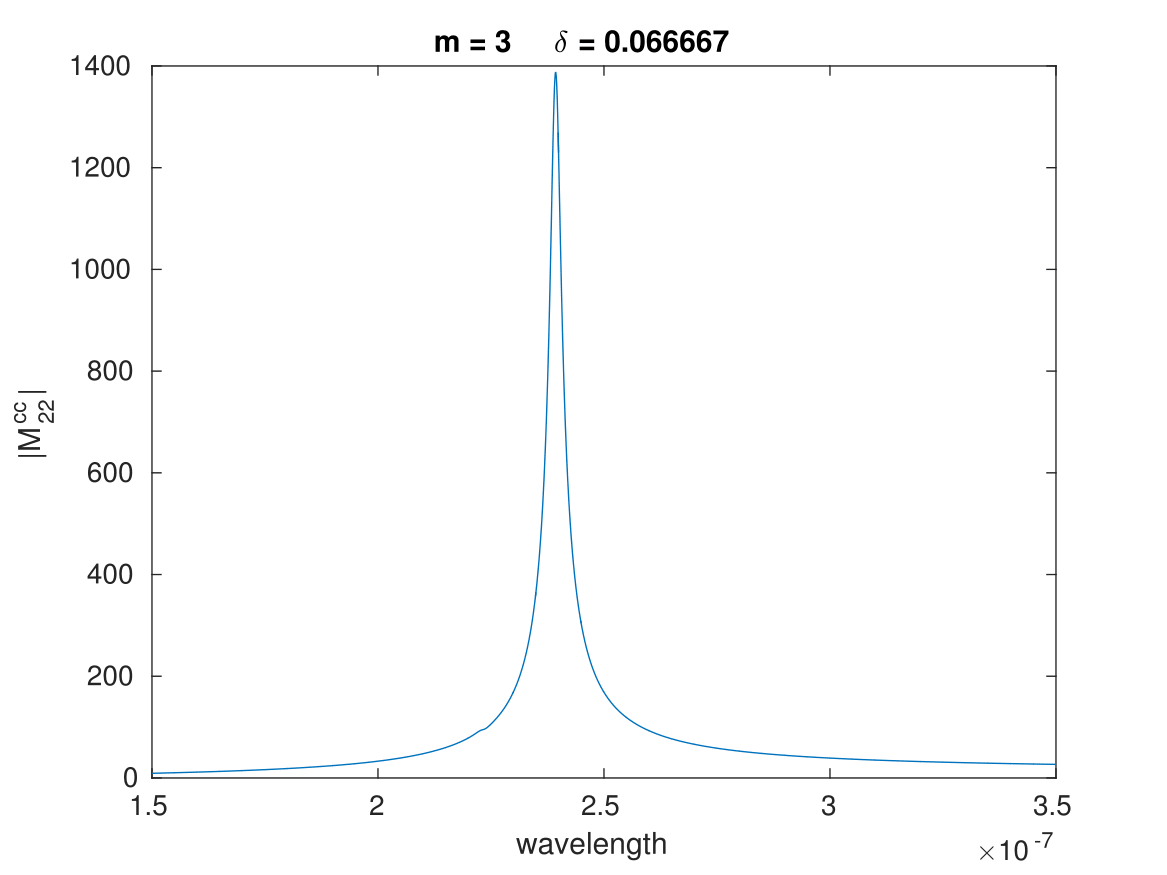}
		\caption{ \label{fig-algbm=3 algebraic} From top to bottom and left to right: initial shape, reconstructed shape, $|m_{11}|$ and $|M_{22}^{cc}|$ with respect to the wavelength for $m = 3$ and $\delta= 0.066667$.}
	\end{center}
\end{figure}

\begin{figure}[!h]
	\begin{center}
		\includegraphics[scale=0.15]{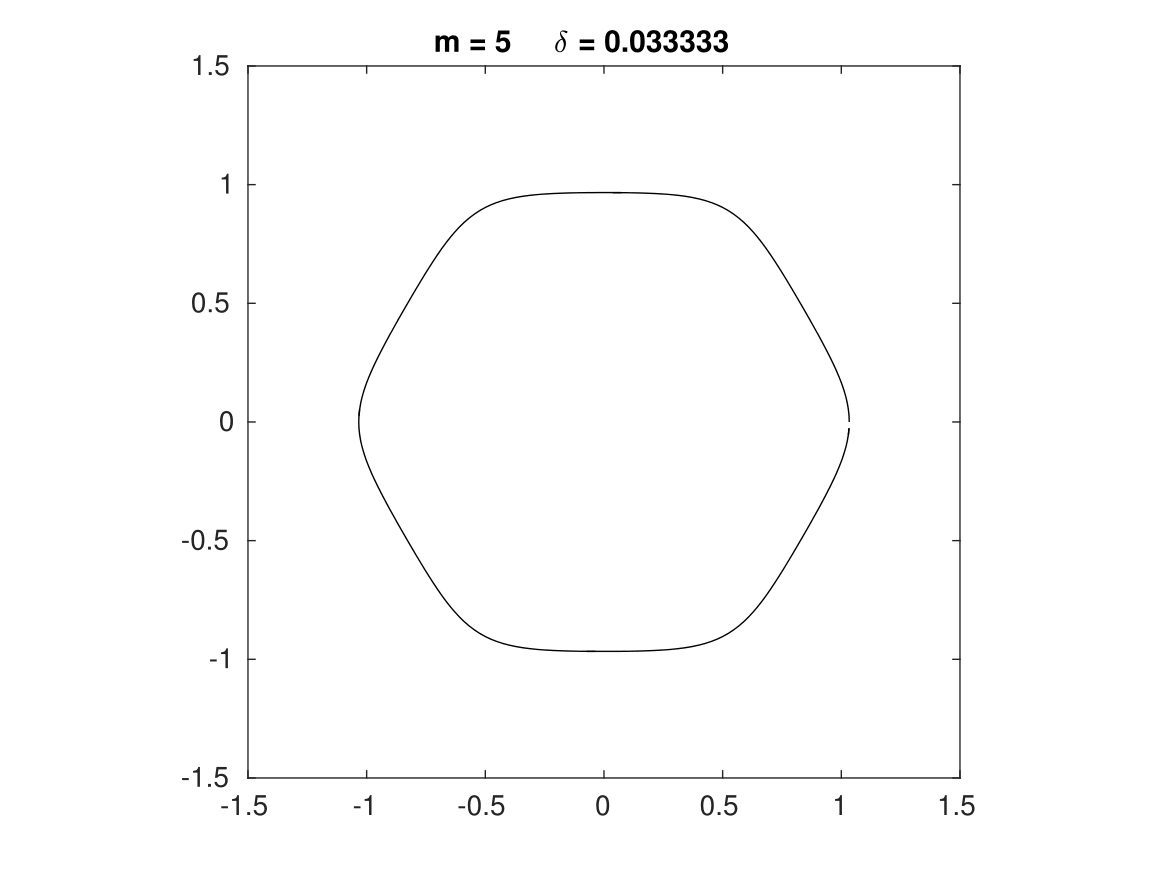}
		\includegraphics[scale=0.15]{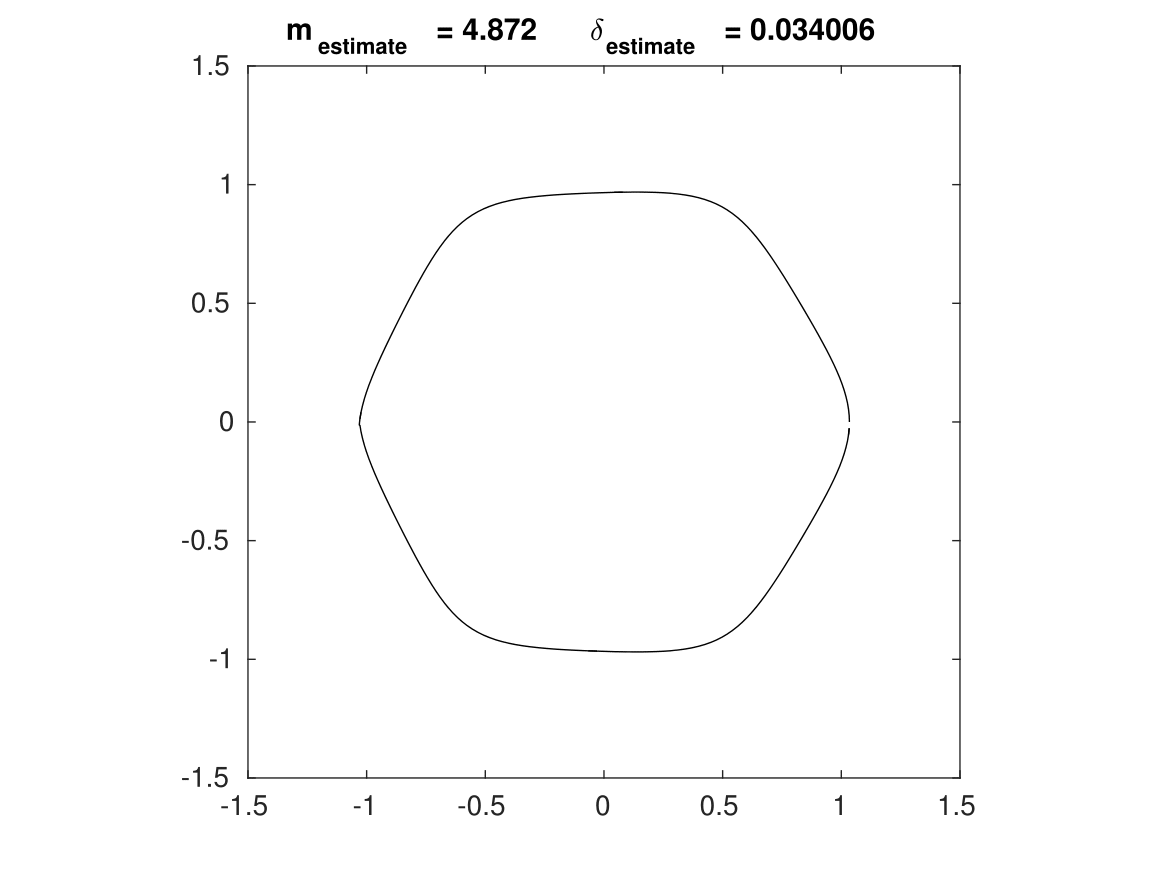}
		\includegraphics[scale=0.15]{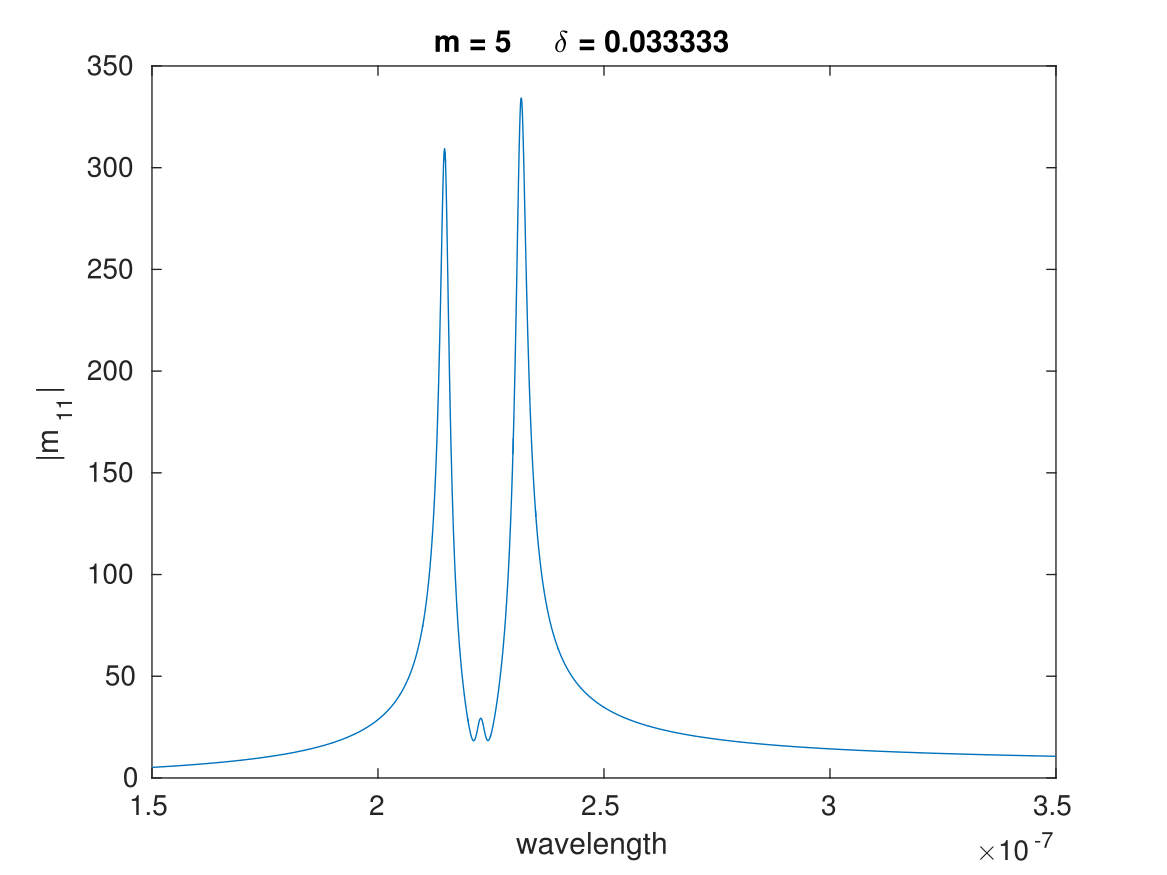}
		\includegraphics[scale=0.15]{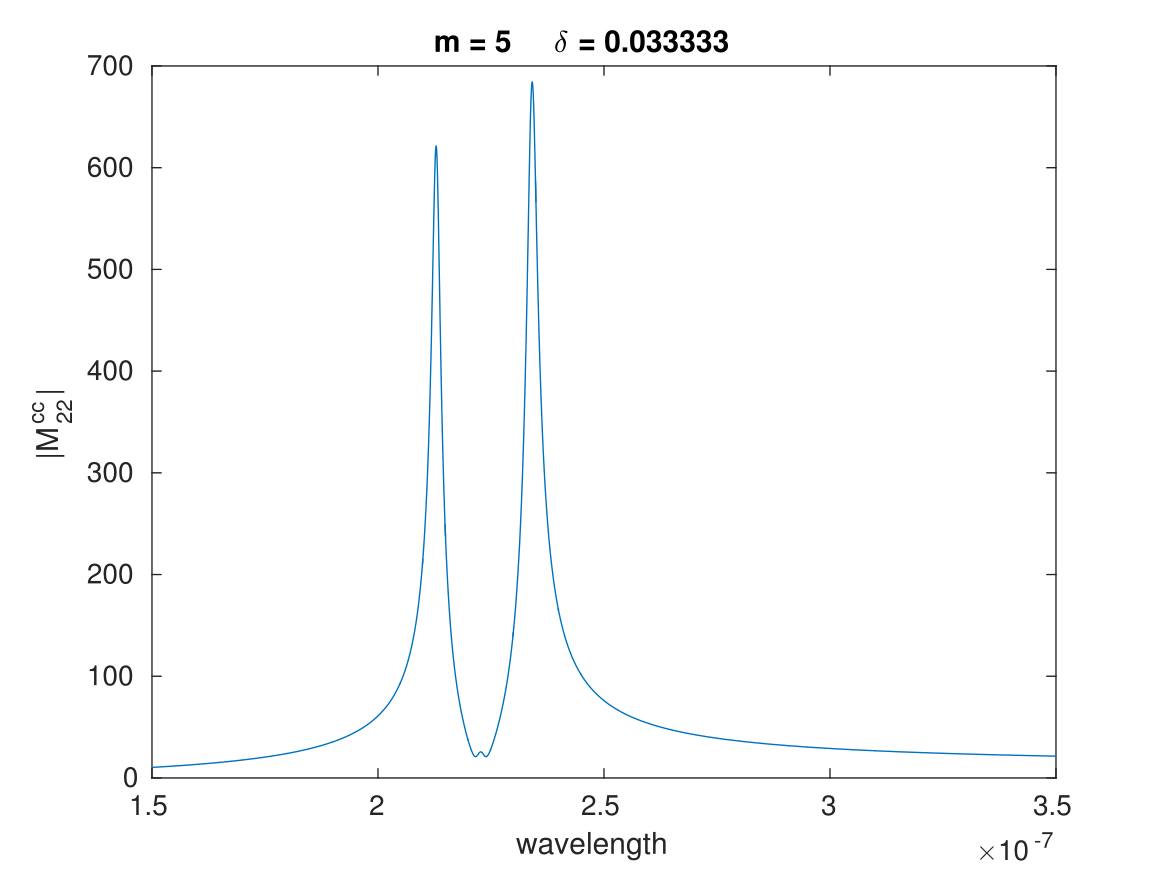}
		\caption{ \label{fig-algbm=5 algebraic} From top to bottom and left to right: initial shape, reconstructed shape, $|m_{11}|$ and $|M_{22}^{cc}|$ with respect to the wavelength for $m = 5$ and $\delta =0.03333$.}
	\end{center}
\end{figure}

\begin{figure}[!h]
	\begin{center}
		\includegraphics[scale=0.15]{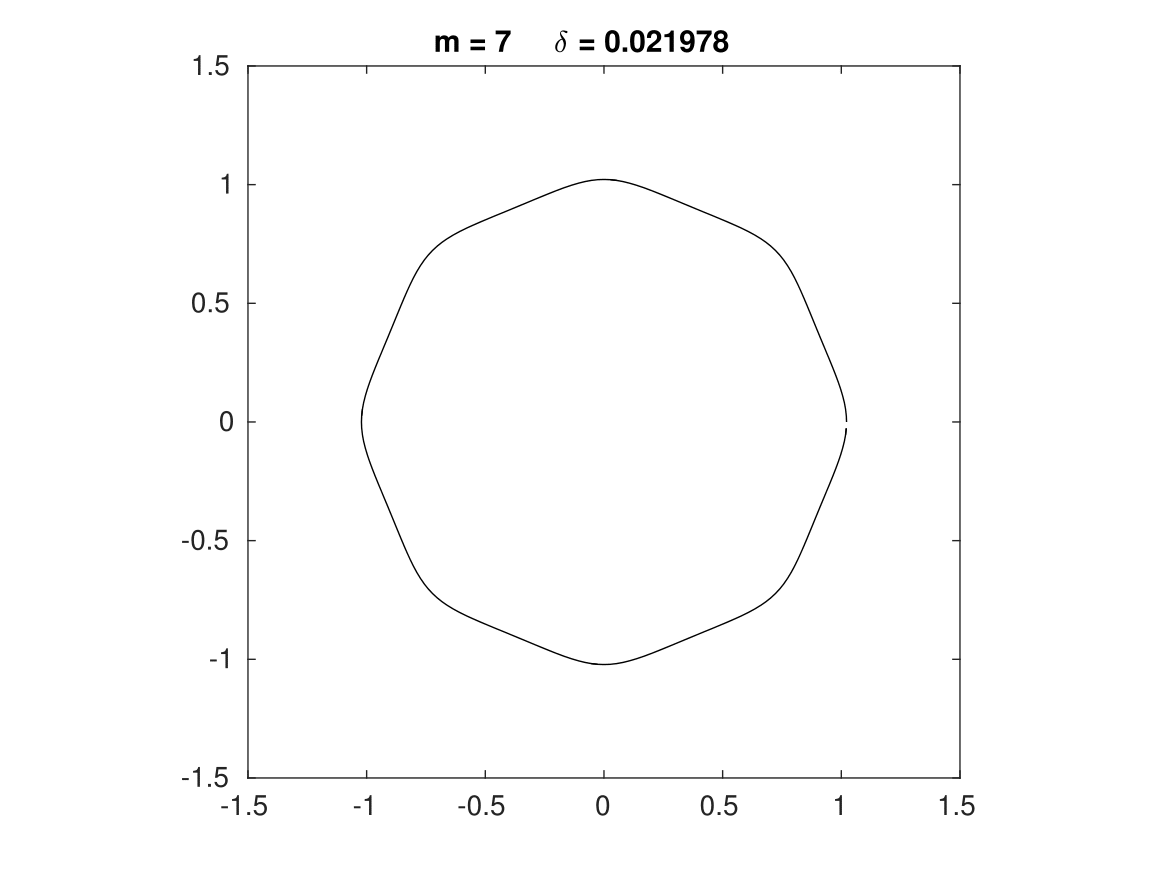}
		\includegraphics[scale=0.15]{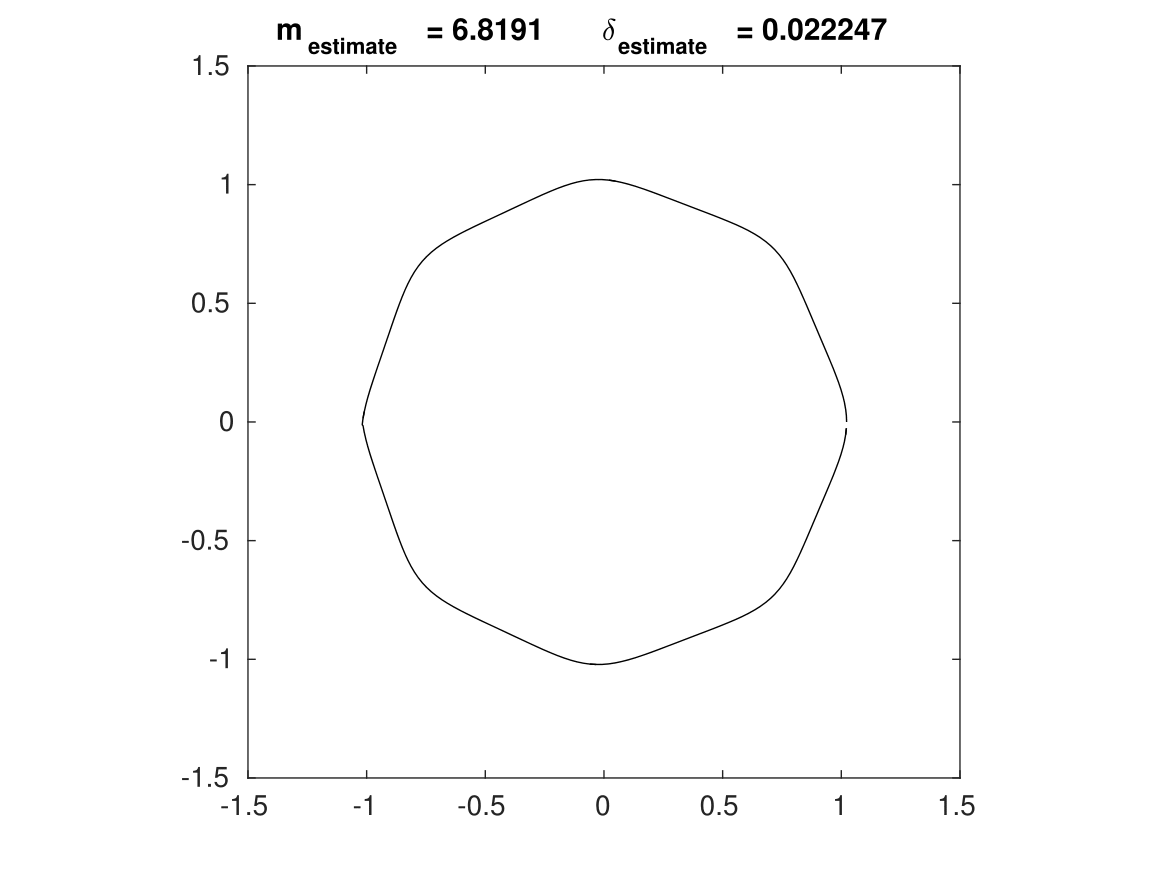}
		\includegraphics[scale=0.15]{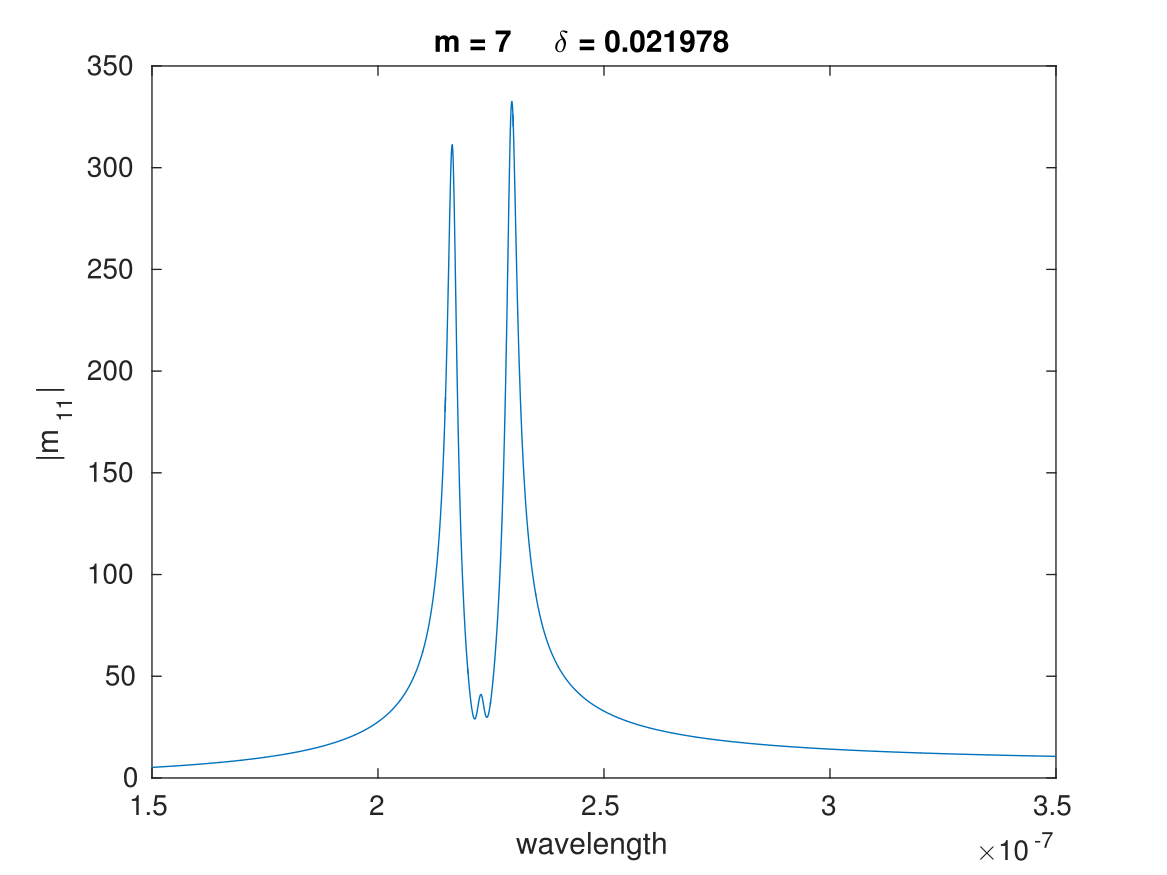}
		\includegraphics[scale=0.15]{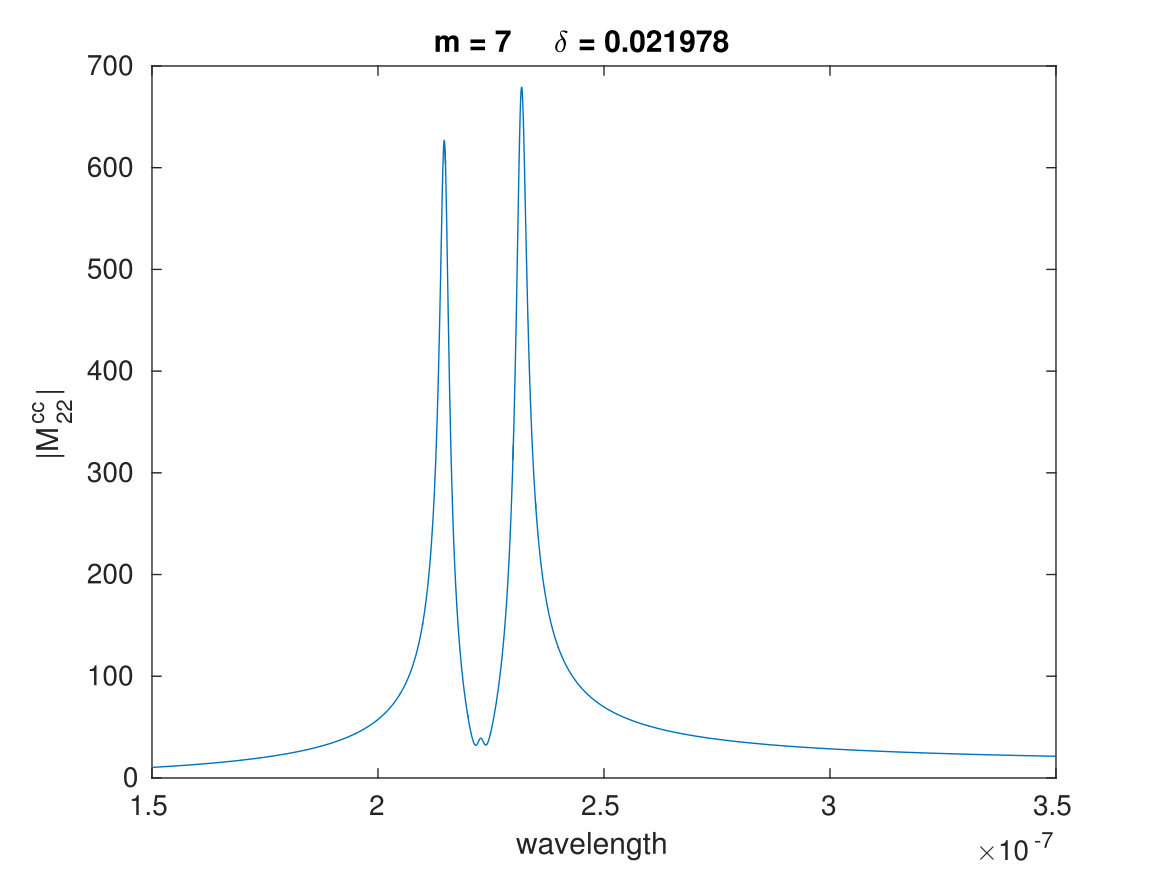}
		\caption{ \label{fig-algbm=7 algebraic} From top to bottom and left to right: initial shape, reconstructed shape, $|m_{11}|$ and $|M_{22}^{cc}|$ with respect to the wavelength for $m = 7$, $\delta = 0.021978$.}
	\end{center}
\end{figure}

Figures \ref{fig-algbm=3 algebraic}, \ref{fig-algbm=5 algebraic}, and \ref{fig-algbm=7 algebraic} present examples of algebraic domains and their reconstructions, where a circle of radius one has been transformed for different values of $m$ and $\delta$.

\begin{figure}[!h]
	\begin{center}
		\includegraphics[scale=0.15]{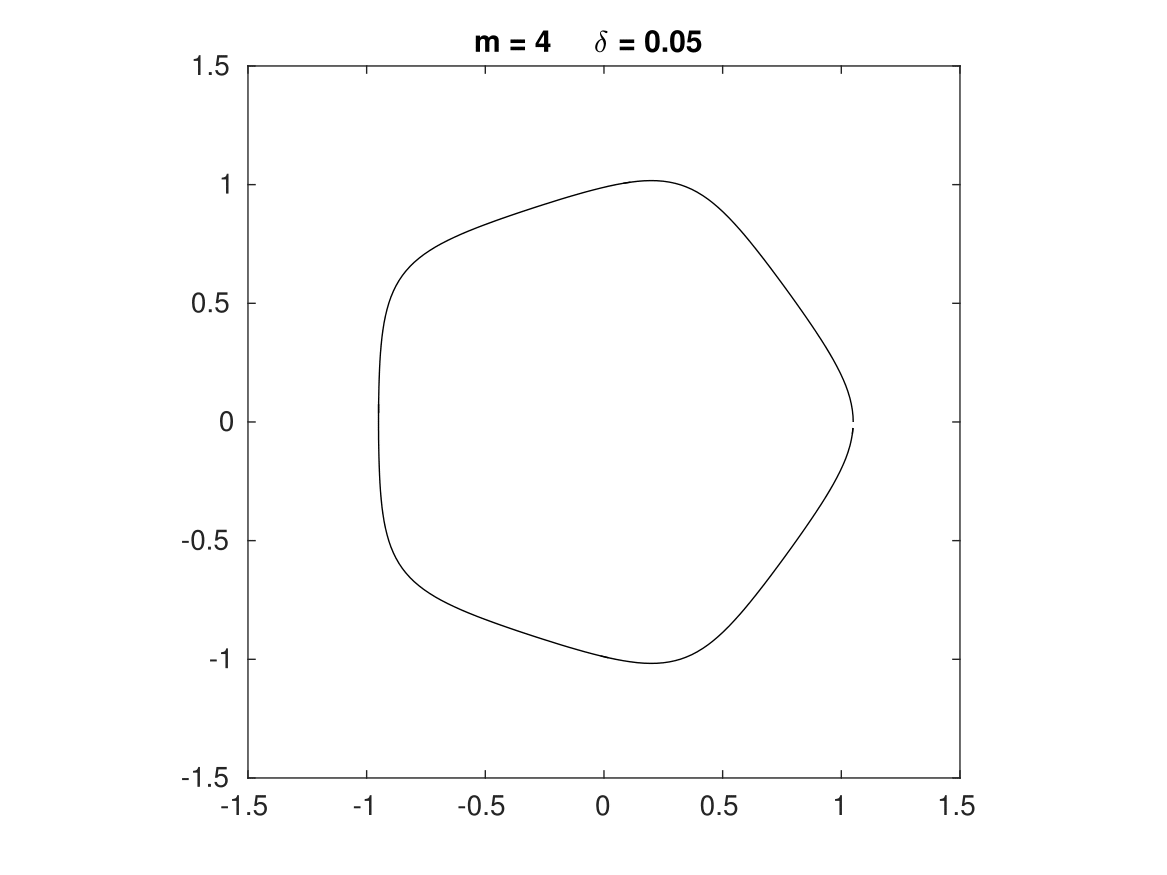}
		\includegraphics[scale=0.15]{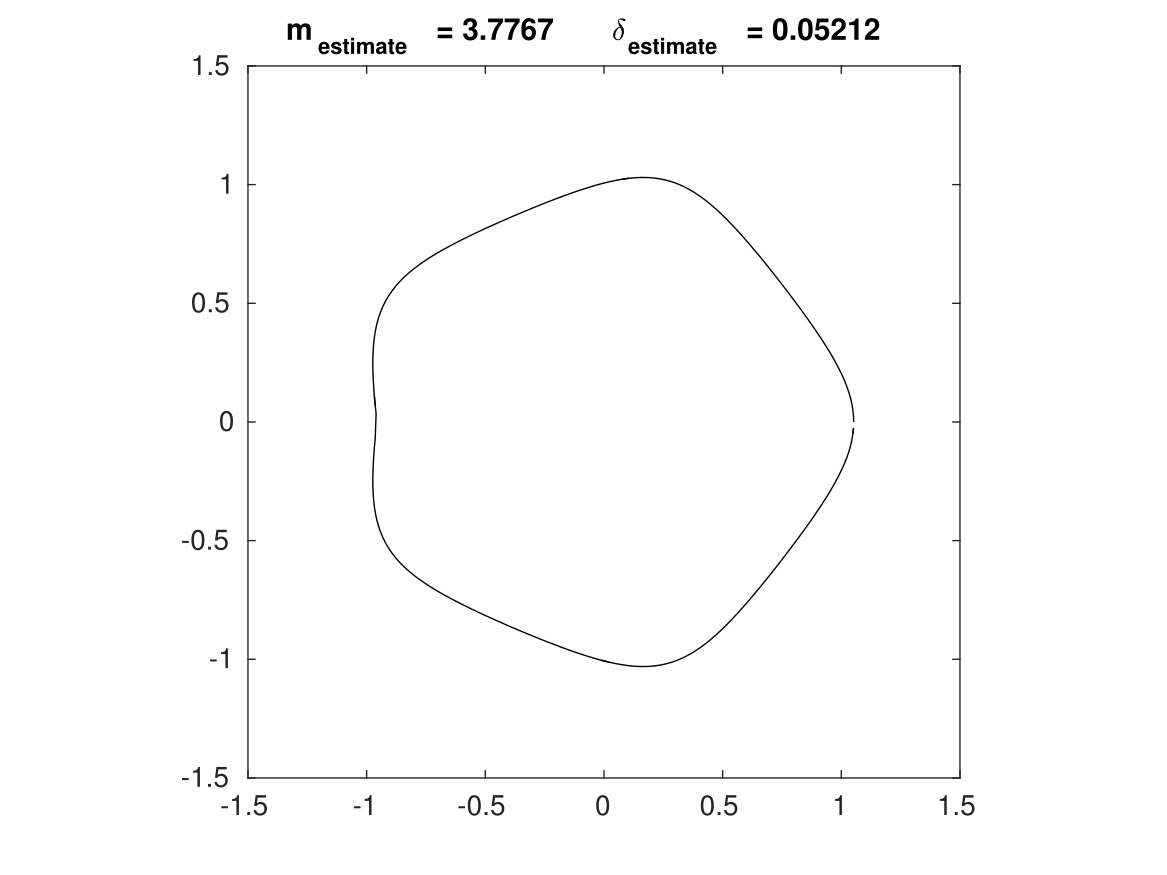}
		\includegraphics[scale=0.15]{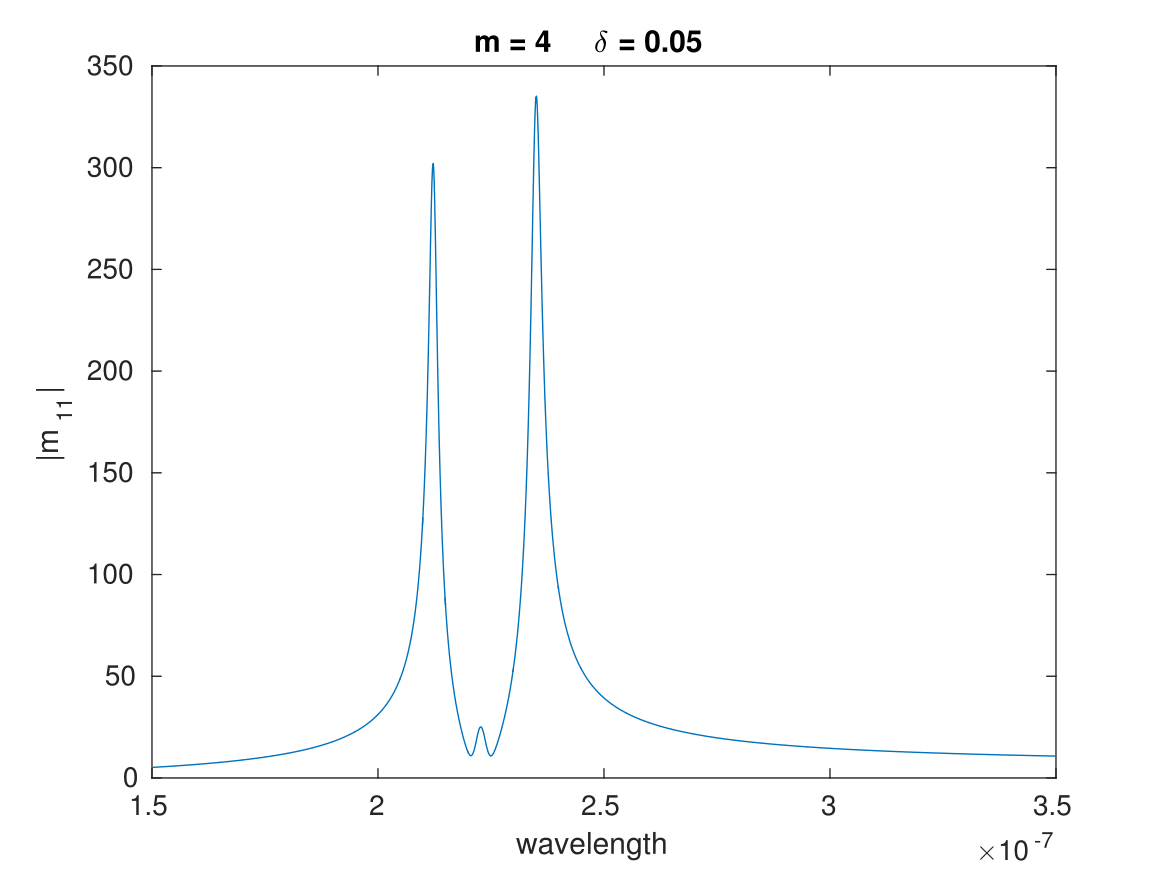}
		\includegraphics[scale=0.15]{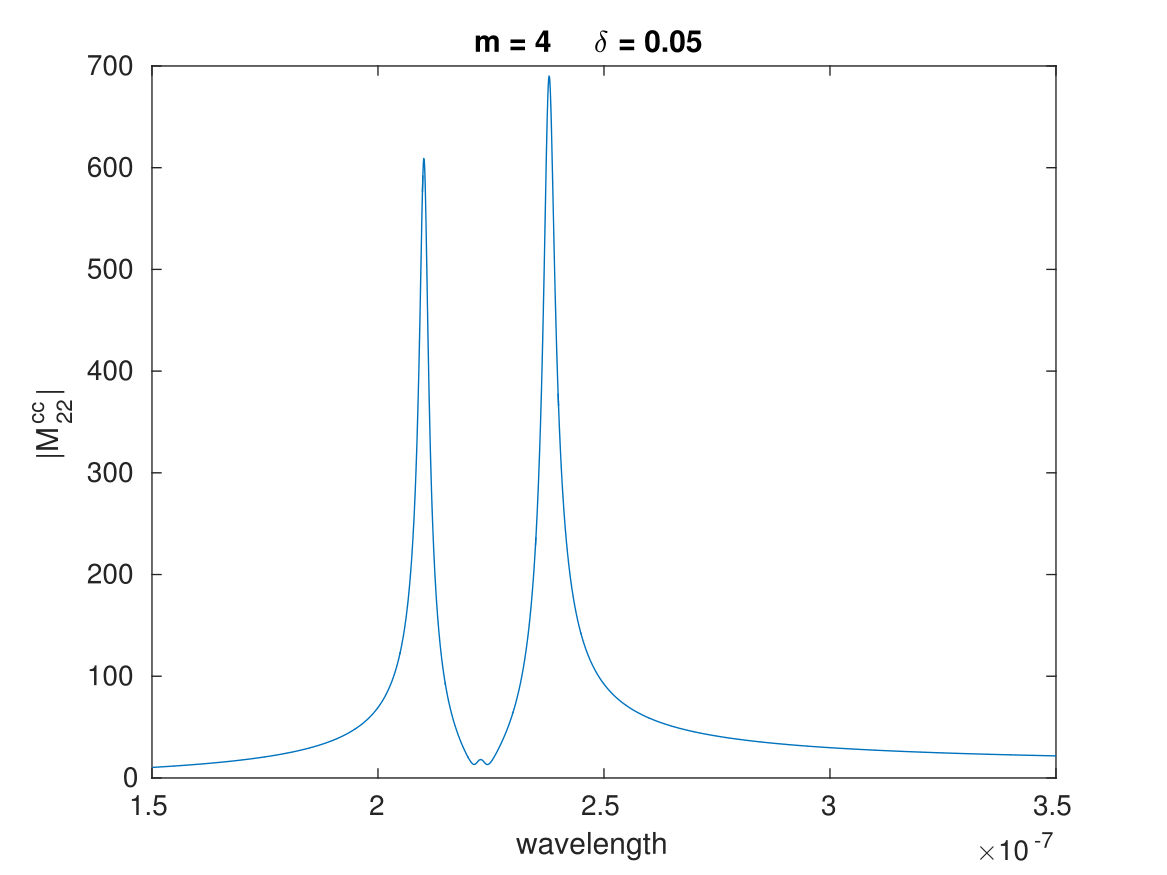}
		\caption{ \label{fig-algbm=4 algebraic} From top to bottom and left to right: initial shape, reconstructed shape, $|m_{11}|$ and $|M_{22}^{cc}|$ with respect to the wavelength for $m = 4$, $\delta = 0.05$.}
	\end{center}
\end{figure}

\begin{figure}[!h]
	\begin{center}
		\includegraphics[scale=0.15]{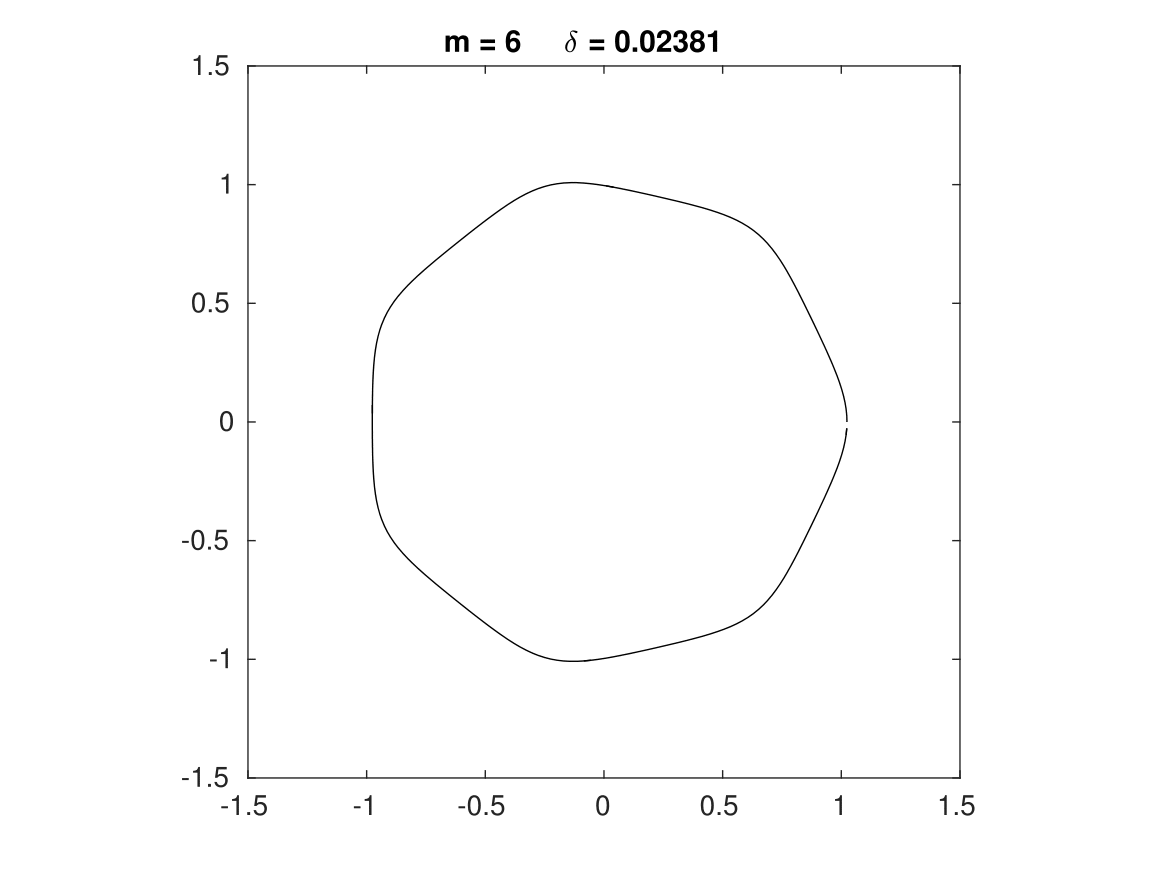}
		\includegraphics[scale=0.15]{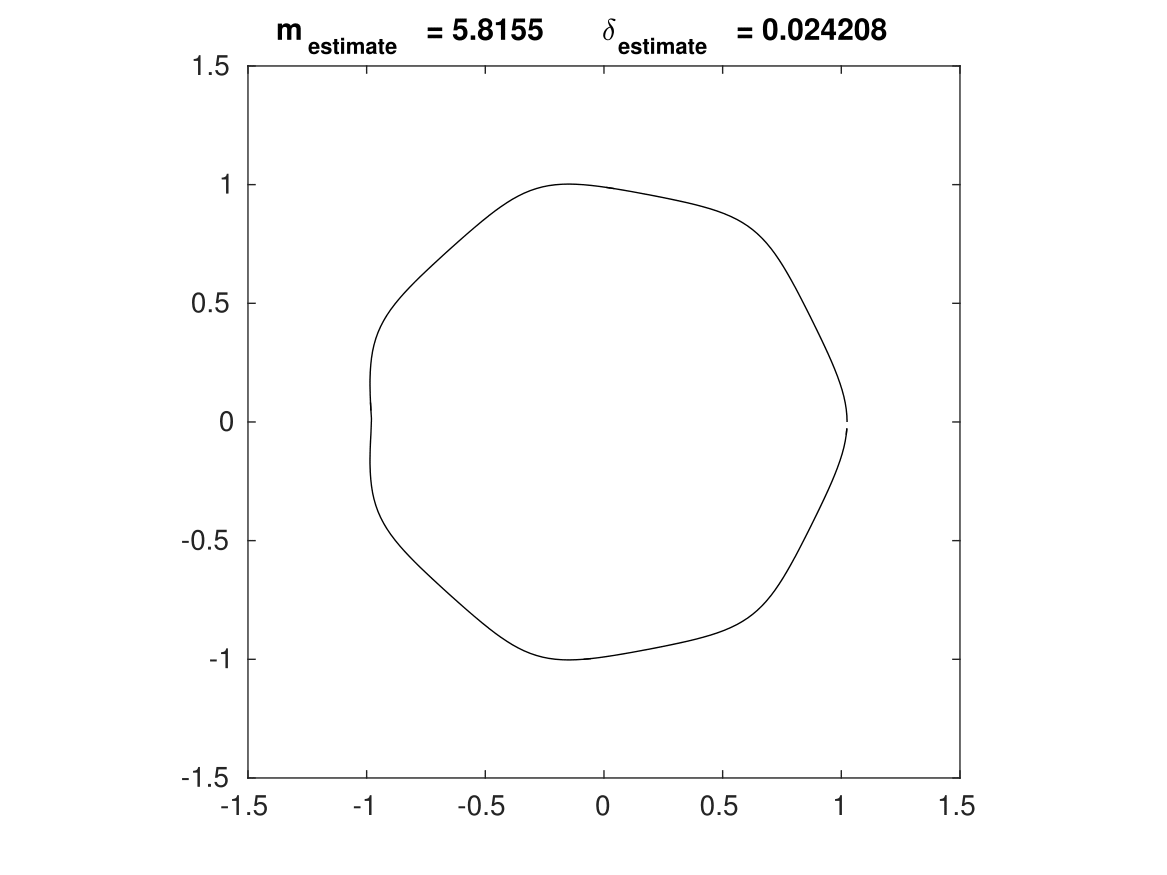}
		\includegraphics[scale=0.15]{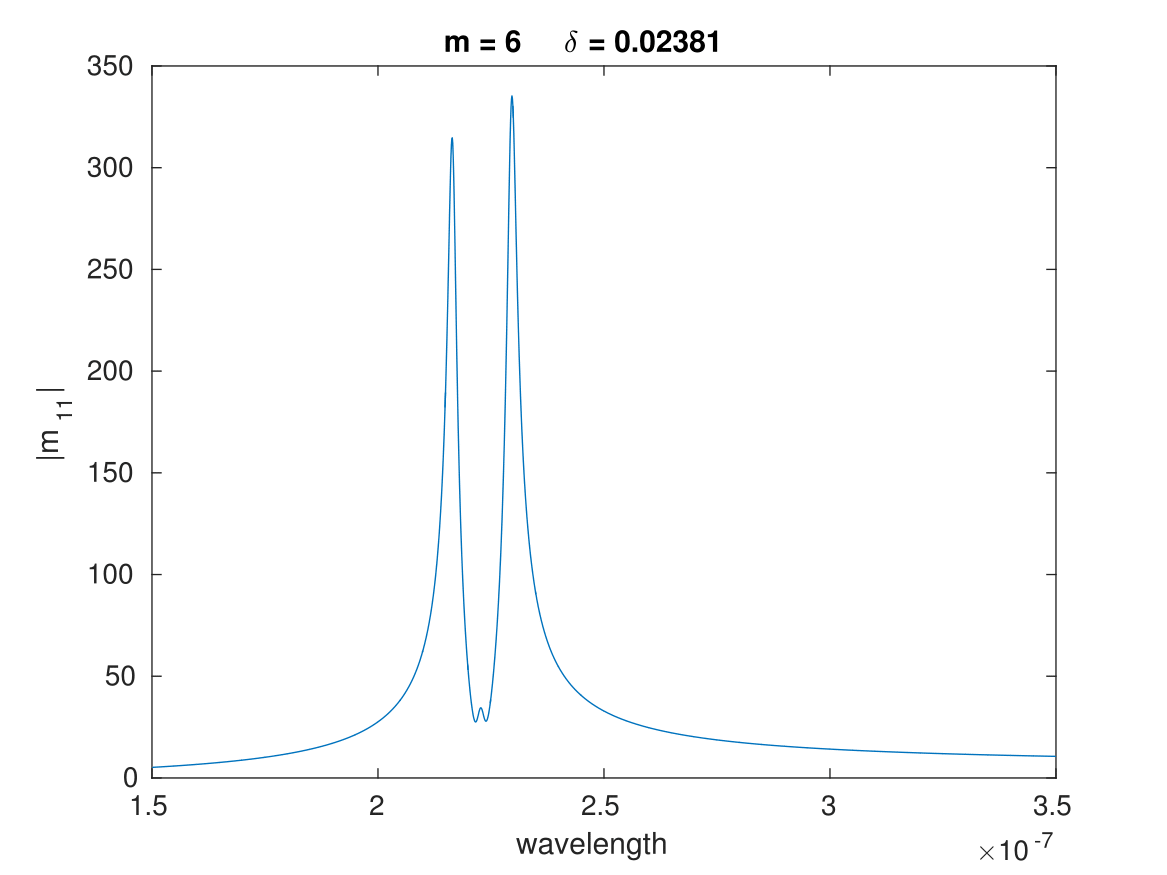}
		\includegraphics[scale=0.15]{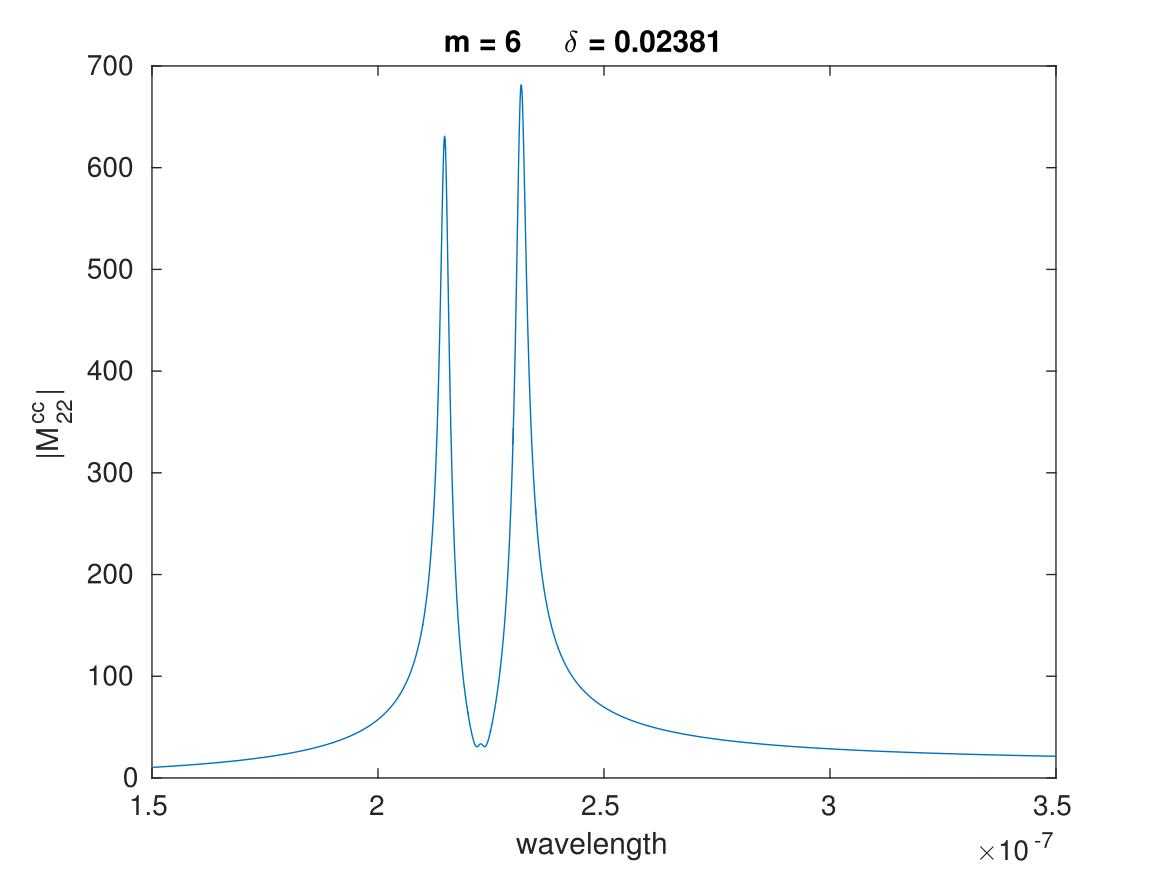}
		\caption{ \label{fig-algbm=6 algebraic} From top to bottom and left to right: initial shape, reconstructed shape, $|m_{11}|$ and $|M_{22}^{cc}|$ with respect to the wavelength for $m = 6$, $\delta = 0.02381$.}
	\end{center}
\end{figure}

Figures \ref{fig-algbm=4 algebraic} and \ref{fig-algbm=6 algebraic} present examples of algebraic domains and their reconstructions, where a circle of radius one has been transformed for $m=4$ and $m=6$ and $\delta=0.02$.

\begin{figure}[h!]
	\begin{center}
		\includegraphics[scale=0.12]{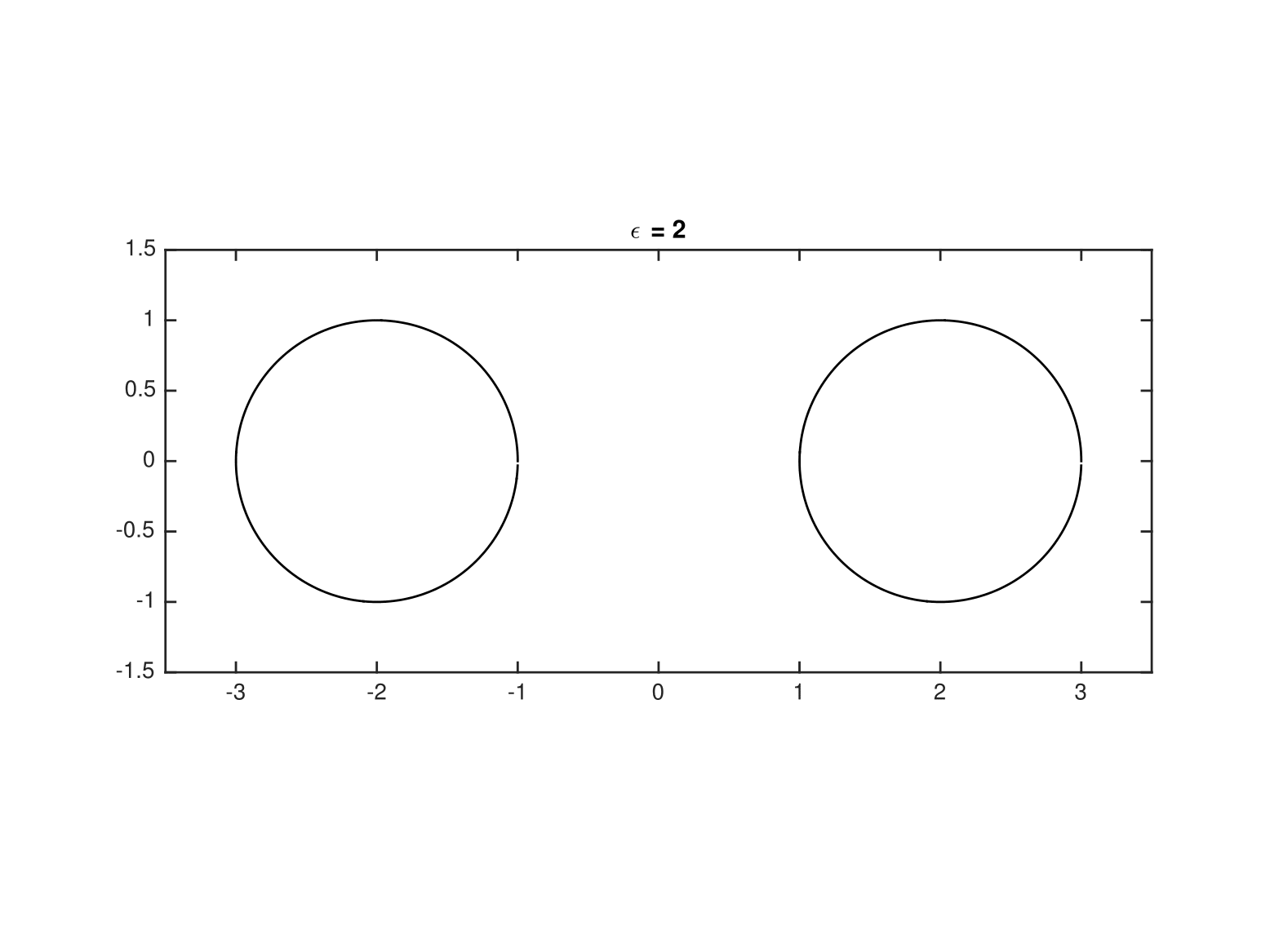}
		\includegraphics[scale=0.12]{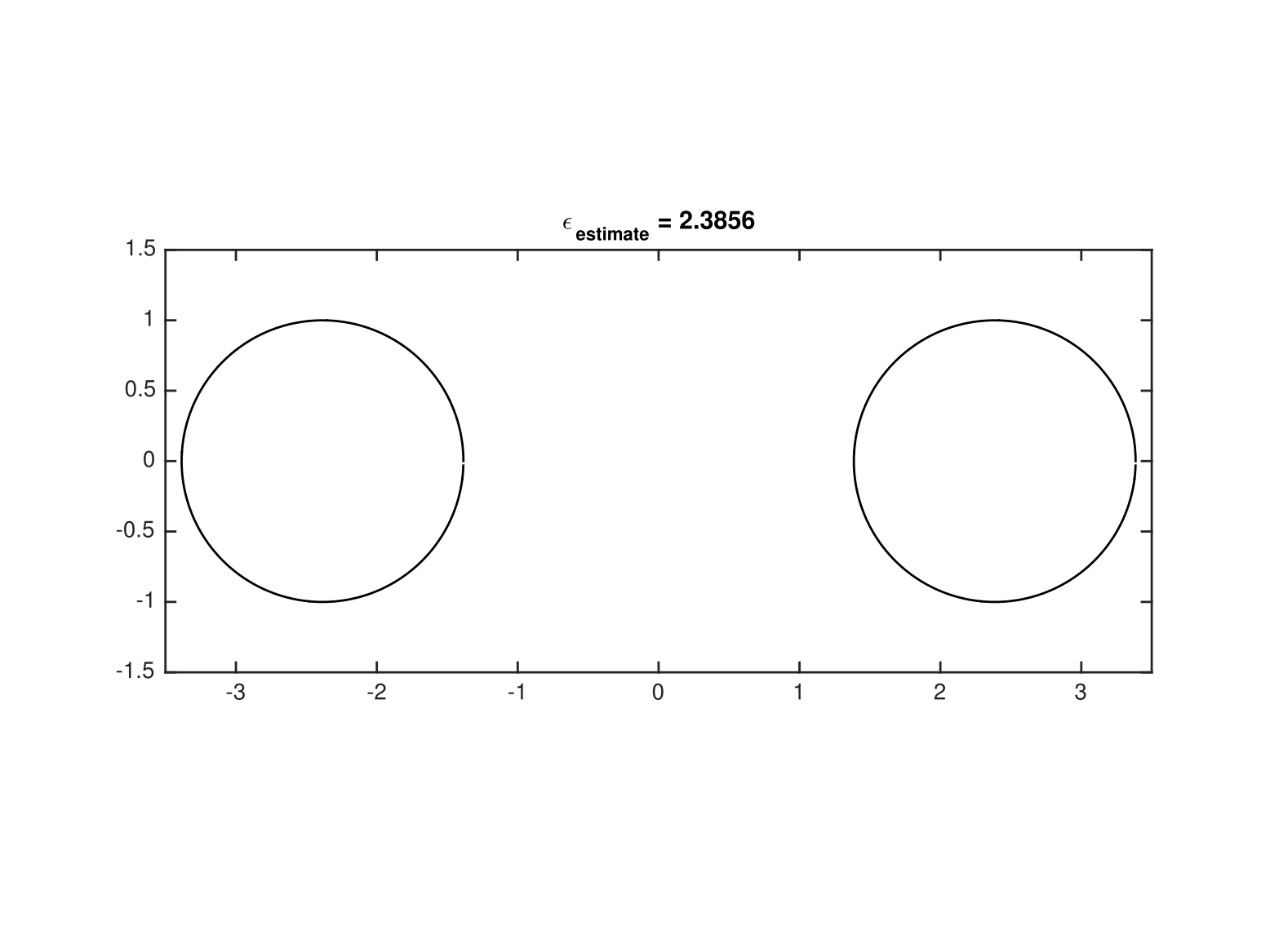}
		\includegraphics[scale=0.12]{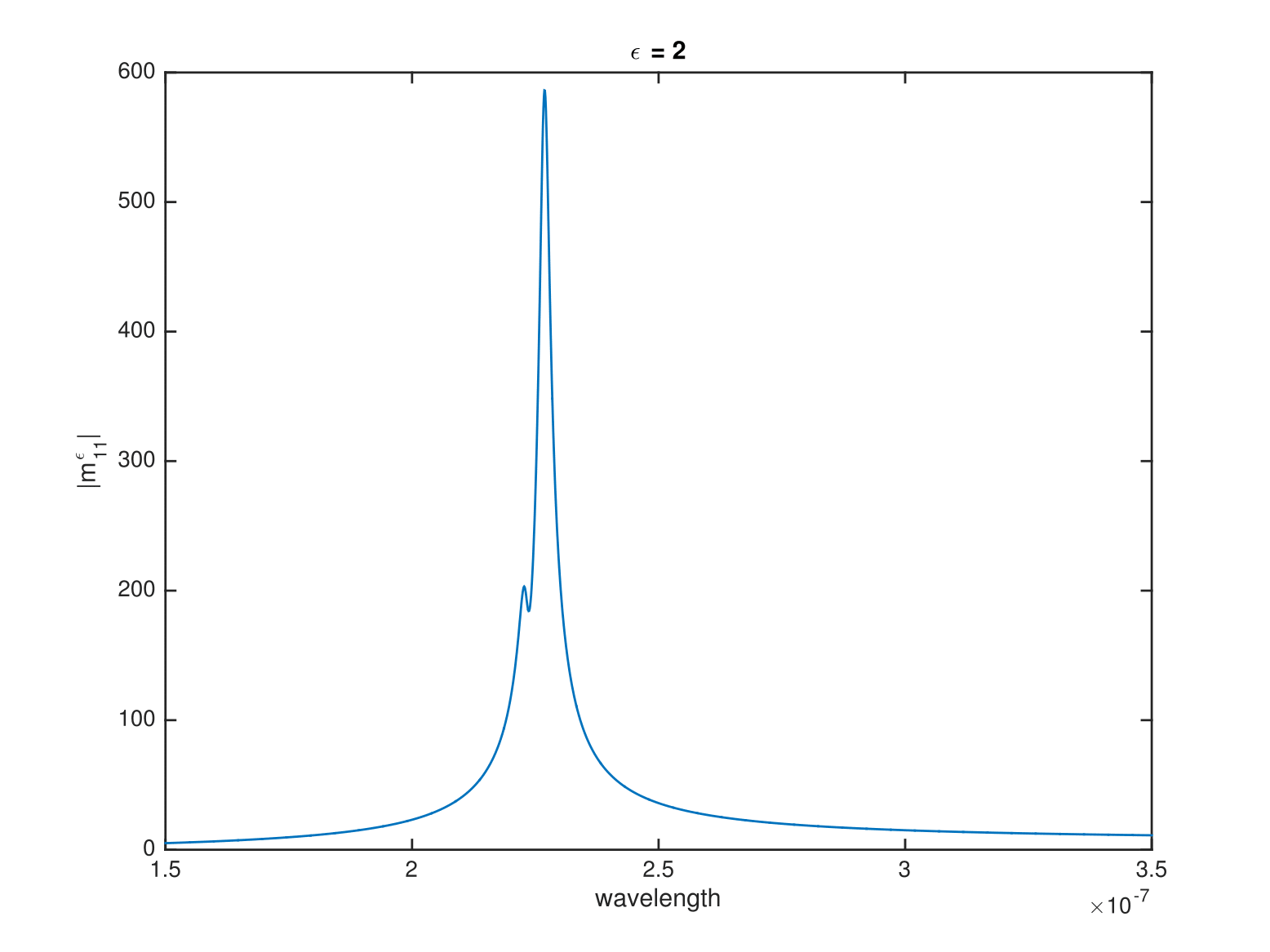}
		\caption{ \label{fig-twocircles2 algebraic} From top left to bottom: initial shape, reconstructed shape and $|m_{11}|$ with respect to the wavelength for $\epsilon = 2$.}
	\end{center}
\end{figure}

\begin{figure}[h!]
 	\begin{center}
 		\includegraphics[scale=0.12]{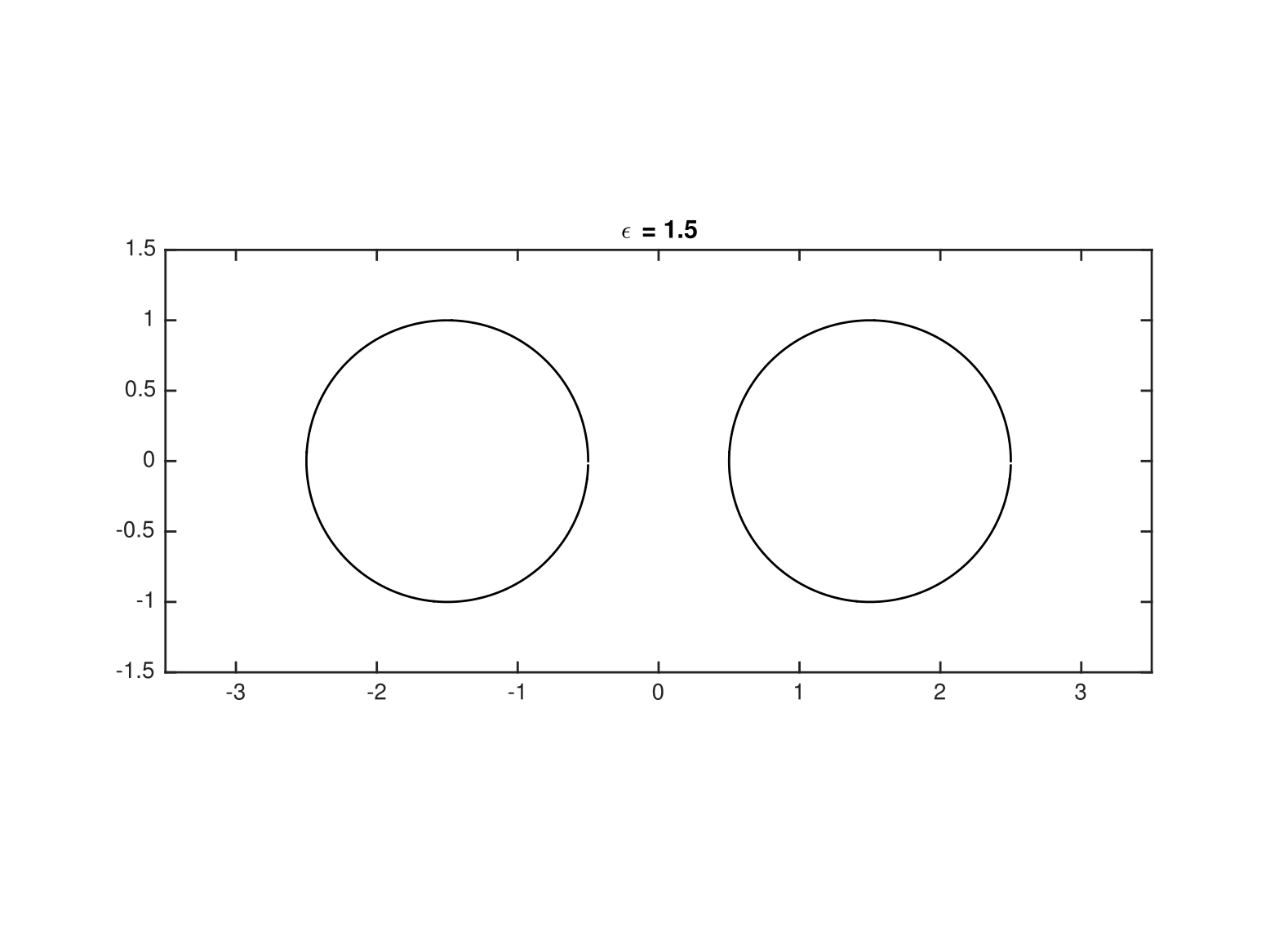}
 		\includegraphics[scale=0.12]{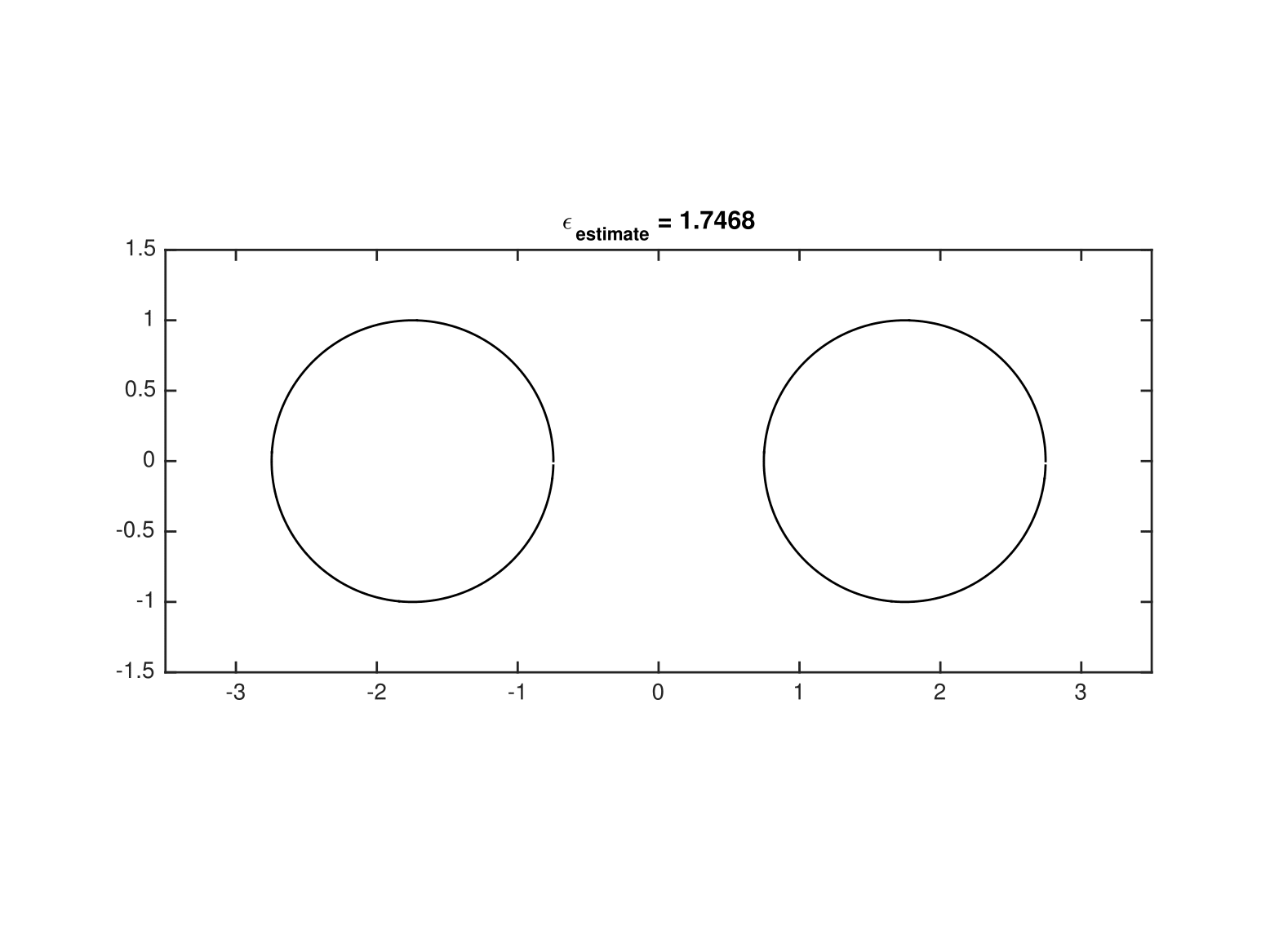}
 		\includegraphics[scale=0.12]{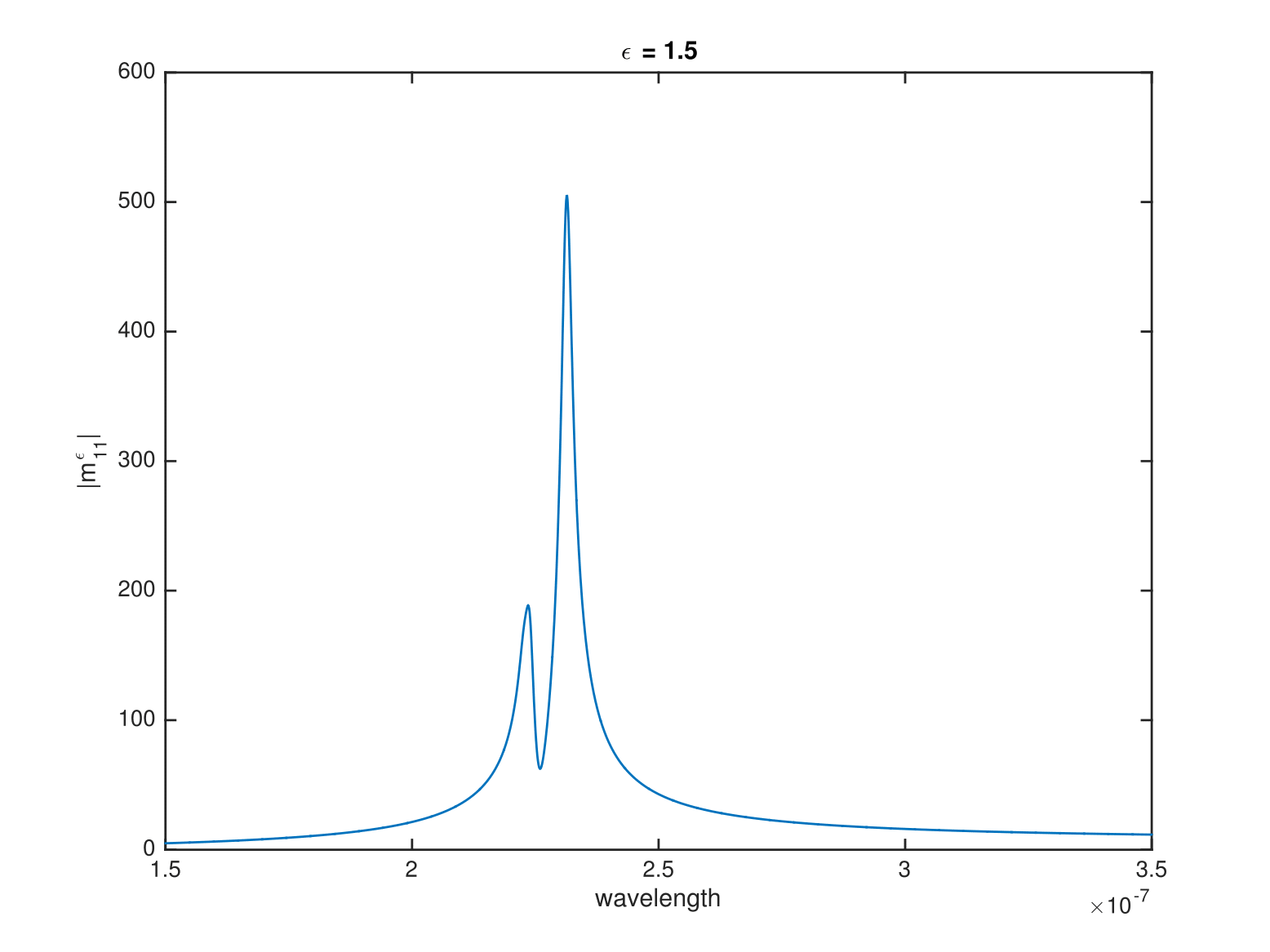}
 		\caption{ \label{fig-twocircles1.5 algebraic} From top left to bottom: initial shape, reconstructed shape and $|m_{11}|$ with respect to the wavelength for $\epsilon = 1.5$.}
 	\end{center}
 \end{figure}

\begin{figure}[h!]
	\begin{center}
		\includegraphics[scale=0.12]{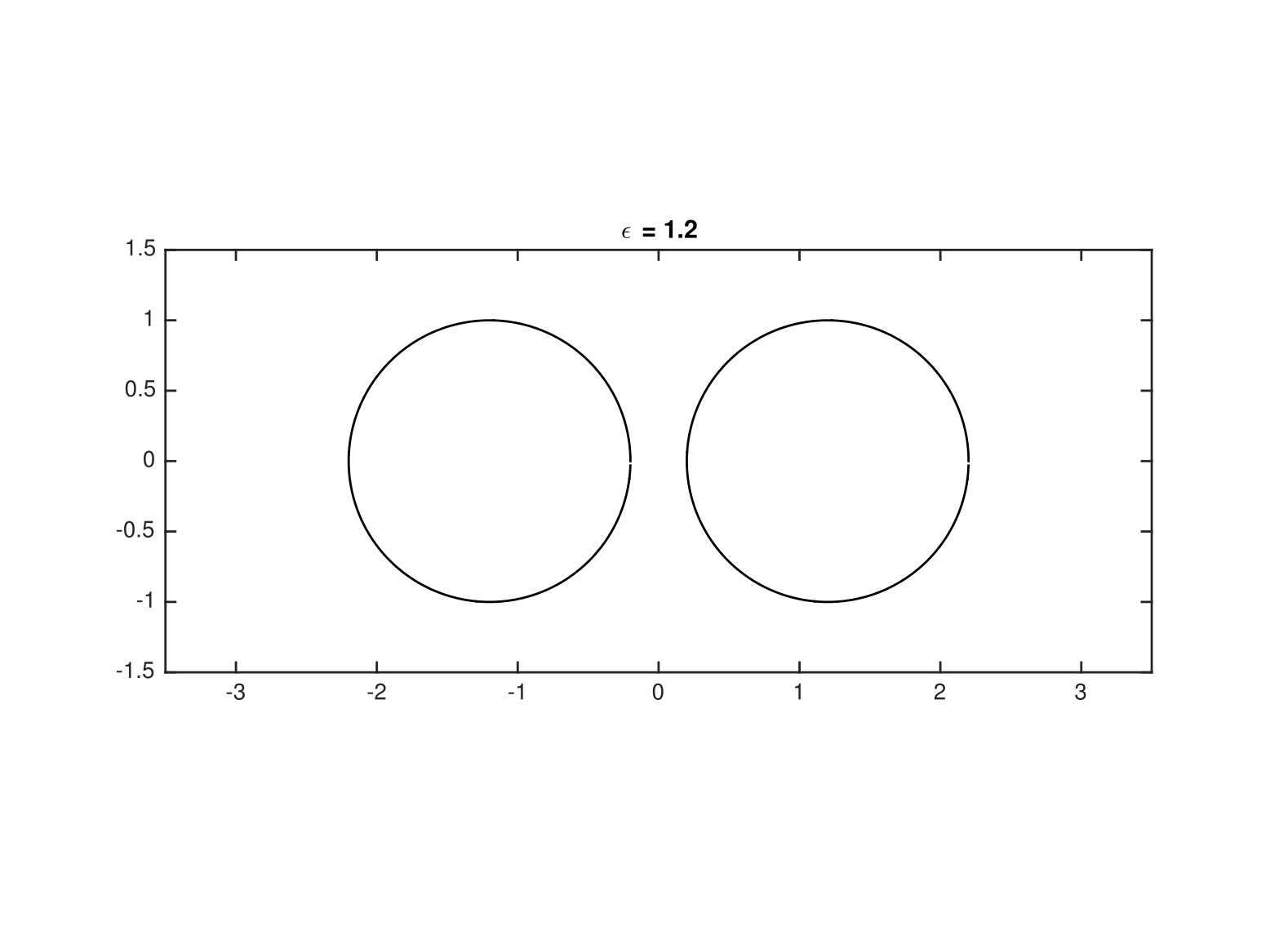}
		\includegraphics[scale=0.12]{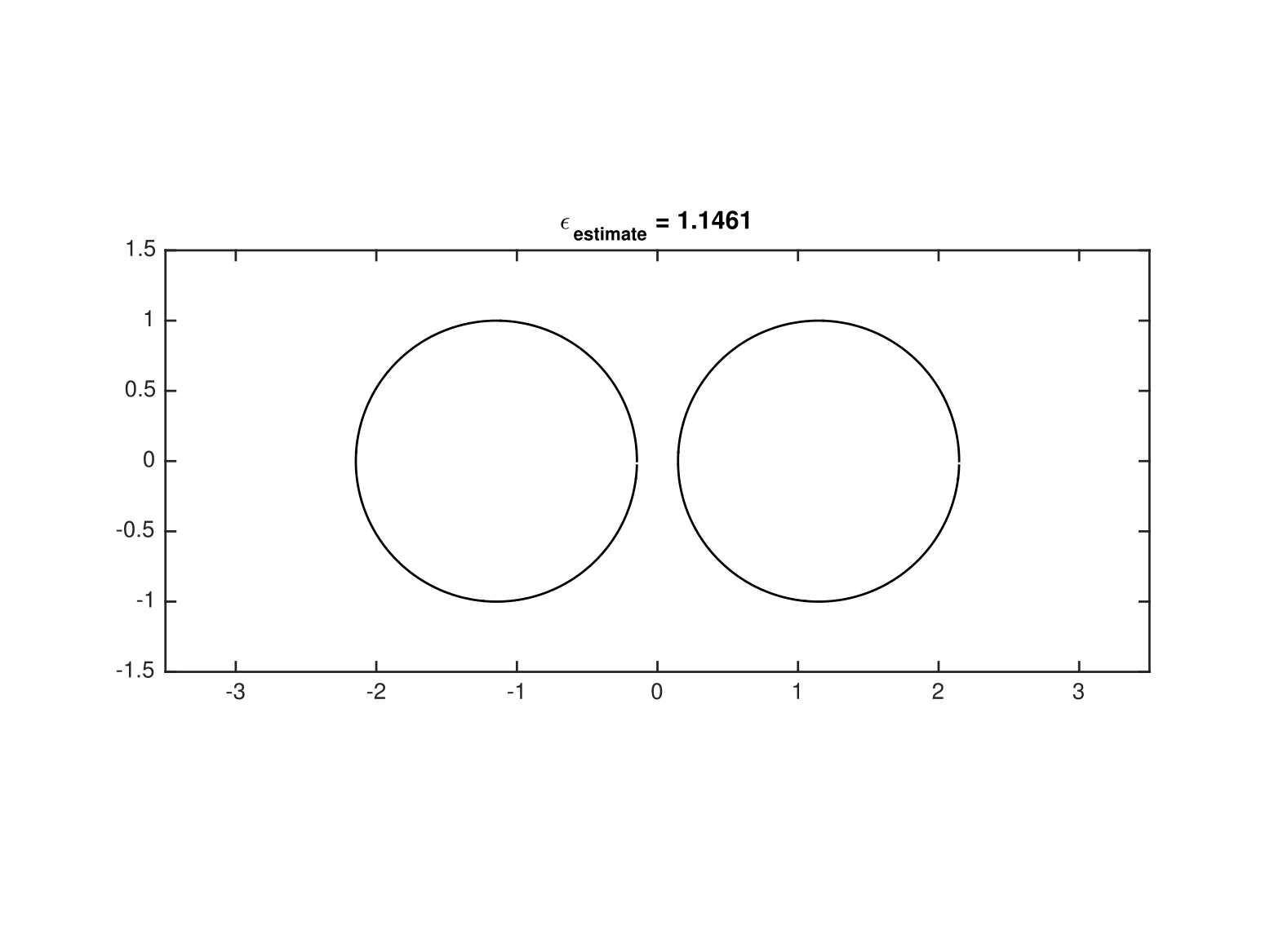}
		\includegraphics[scale=0.12]{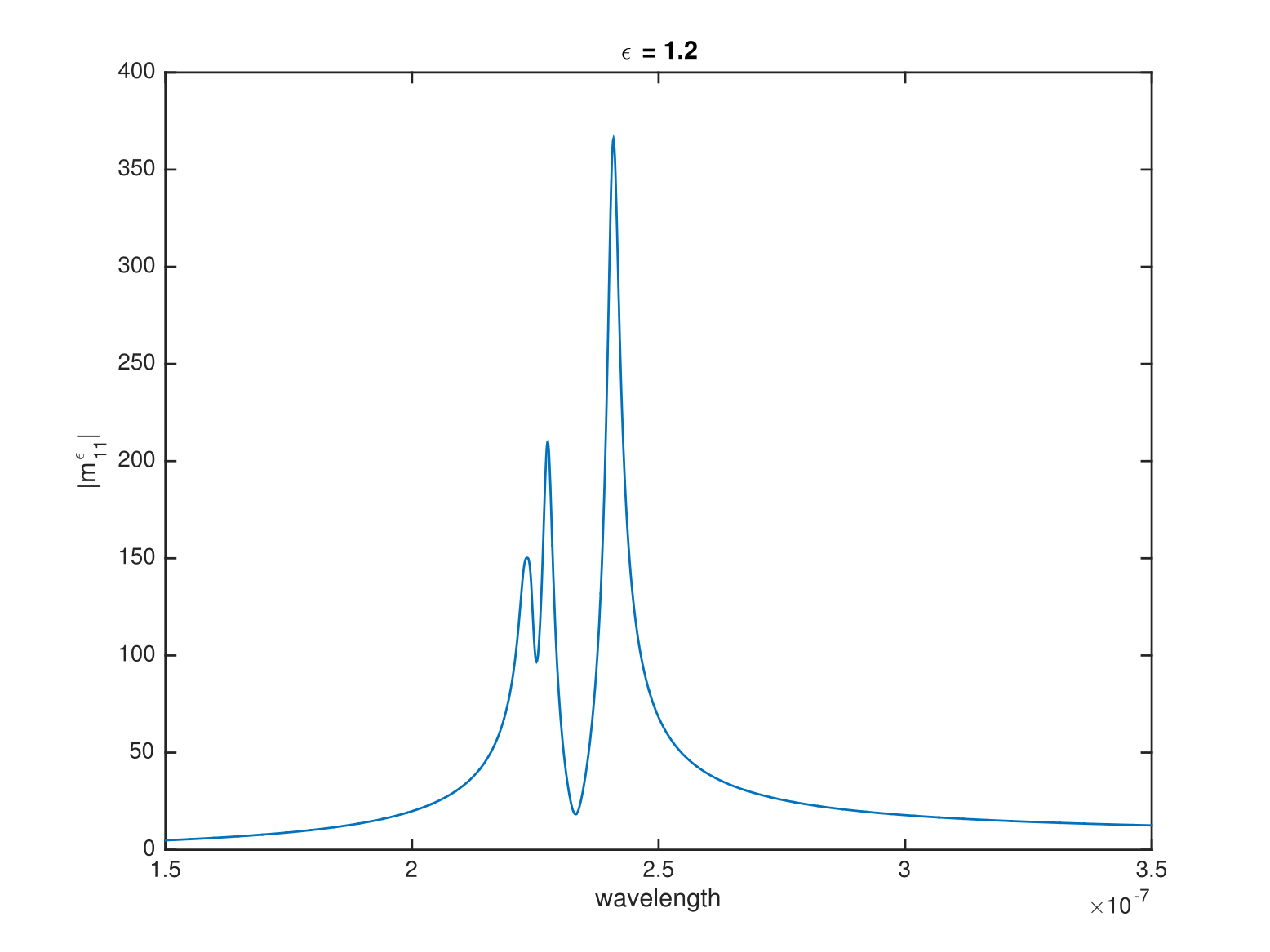}
		\caption{ \label{fig-twocircles1.2 algebraic} From top left to bottom right: initial shape, reconstructed shape and $|m_{11}|$ with  respect to the wavelength for $\epsilon = 1.2$.}
	\end{center}
\end{figure}

Figures \ref{fig-twocircles2 algebraic}, \ref{fig-twocircles1.5 algebraic}, and \ref{fig-twocircles1.2 algebraic} show examples of two circles of radius one separated by a distance $\ep$, and their reconstructions.

\section{Concluding remarks}
In this paper we have proved for a class of algebraic domains that the associated plasmonic resonances 
can be used to classify them. It would be very interesting to prove a similar result for all quadrature domains or all algebraic domains. We have also reconstructed the separation distance between two nanoparticles of circular shape from measurements of their first collective plasmonic resonances. Another challenging problem would be to generalize this result to more components and arbitrary shaped particles. 

\appendix 

\section{Blow-up of the gradient of the potential at plasmonic resonances for two nearly touching disks}
\label{appendixA}
In this appendix, we consider two nearly touching disks and estimate the blow up of the gradient of the potential for two nearly touching disks at plasmonic resonances. We generalize to the plasmonic case, the estimates derived in \cite{akl, jde} for nearly touching disks with degenerate conductivities. In the nondegenerate case, we refer to \cite{ens,sinum, bonnetier,yli}. 

Let $u$ be the electric potential when an external potential $H$ is applied. In other words, $u$ satisfies
$$
\begin{cases}
\ds\Delta u = 0 \quad& \mbox{ in } B_1 \cup B_2\cup \mathbb{R}^2\setminus \overline{(B_1\cup B_2)},\\
\ds u|_+=u|_- \quad& \mbox{ on } \p B_j,  j=1,2,\\
\ds \frac{\p u}{\p \nu}\Big|_+=k\frac{\p u}{\p \nu}\Big|_- \quad&\mbox{ on } \p B_j, j=1,2,\\
\ds u(\Bx)-H(\Bx)= O(|\Bx|^{-1}) \quad& \mbox{ as } |\Bx| \rightarrow \infty.
\end{cases}
$$
$u$ can be represented as follows:
$$
u=H+\mathbb{S}[\varphi],
$$
where $\varphi$ is the solution to
$$
(\lambda \mathbb{I}-\mathbb{K}^*)[\varphi] = 
\begin{bmatrix}
\ds \p_\nu H|_{\p B_1}
\\
\ds \p_\nu H|_{\p B_2}
\end{bmatrix}.
$$
Suppose that $H(x,y)=E_0 x$.
From \eqnref{x1_bipolar_2}, we obtain
\begin{align}
\frac{\p H}{\p\nu}|_{\xi=\pm\xi_0} &=\pm  E_0 h(\xi_0,\theta)\alpha \sum_{m=-\infty}^\infty (-|m|)e^{-|m|\xi_0+i m\theta}.
\end{align}
On the other hand, by \eqnref{normal_eigenvec},  
\begin{align}
\begin{bmatrix}
\ds \p_\nu H|_{\p B_1}
\\
\ds \p_\nu H|_{\p B_2}
\end{bmatrix}
=\sum_{m=-\infty}^\infty E_0 \alpha |m|{\sqrt{2\pi(1-e^{-2|m|\xi_0})}} e^{-|m|\xi_0} \Psi^+_{\ep,m} +0\cdot \Psi_{\ep,m}^-.
\end{align}
Hence, 
\begin{align}
\varphi=\sum_{n \neq 0} \frac{1}{\lambda-\lambda_{\ep,n}^+} \Big(E_0\alpha |n|{\sqrt{2\pi(1-e^{-2|m|\xi_0})}} e^{-|m|\xi_0}\Big)\Psi_{\ep,n}^+.
\end{align}
Then \eqnref{Single_normal_eigen} yields
\begin{align}
u=E_0 x+ E_0\sum_{n\neq 0} \frac{\alpha e^{-2|n|\xi_0}}{\lambda-\lambda_{\ep,n}^+}  \sinh|n|\xi \, e^{i n \theta}.
\end{align}
Now we compute the electric field at the center of the gap $\nabla u(\xi=0,\theta=\pi)$. We have
$$
\nabla u(0,\pi) =  E_0\mathbf{e}_x + E_p \mathbf{e}_x, \quad E_p = E_0 \sum_{n=1}^\infty \frac{4|n| e^{-2|n|\xi_0}}{\lambda-\lambda_{\ep,n}^+} (-1)^n.
$$
Let $k_{\ep,n}^+=-\coth |n| \xi_0$. Note that 
$$
\lambda_{\ep,n}^+ = \frac{k_{\ep,n}^+ +1}{2(k_{\ep,n}^+ -1)}.
$$
$E_p$ can be rewritten as
$$
E_p = E_0\sum_{n =1}^\infty  \frac{(1-k)(k_{\ep,n}^+-1)}{k-k_{\ep,n}^+} 4|n|e^{-2|n|\xi_0}(-1)^n.
$$
Let us assume that $k$ is given by
$$
k = k_{\ep,N}^+ + i\delta, 
$$
for some $N\in\mathbb{N}$, where $\delta>0$ is a small parameter. This implies that, if $\delta$ goes to zero, then the plasmon resonance occurs at the $N$-th mode. For small $\delta>0$, $E_p$ can be approximated by
\begin{align}
E_p \approx E_0\frac{(1-k)(k_{\ep,N}^+-1)}{k-k_{\ep,N}^+} 4N e^{-2N\xi_0}(-1)^N.
\end{align}
Since $k_{\ep,N}^+ \approx - \frac{\sqrt{r}}{N \sqrt{\ep}}$ for small $\ep>0$, we have
\begin{align}
E_p \approx i E_0 \frac{r}{N \delta \ep}4 e^{-2N\xi_0}(-1)^N.
\end{align}
It follows that
$$
\nabla u(0,\pi) \approx iE_0\frac{r}{N \delta \ep}4 e^{-2N\xi_0}(-1)^N \mathbf{e}_x.
$$

\end{document}